\documentclass[12pt]{amsart}

\usepackage{amssymb,amsmath}
\usepackage[all]{xy}

\usepackage{amsmath}
\usepackage{amssymb}
\usepackage{url}
\usepackage{amscd}
\usepackage{mathrsfs}
\usepackage[dvips]{graphicx}

\usepackage[OT2,T1]{fontenc}
\DeclareSymbolFont{cyrletters}{OT2}{wncyr}{m}{n}
\DeclareMathSymbol{\Sha}{\mathalpha}{cyrletters}{"58}

\theoremstyle{plain}
\newtheorem{thm}{Theorem}

\newtheorem{cor}{Corollary}
\newtheorem{lemma}{Lemma}
\newtheorem{proposition}[thm]{Proposition}

\newcommand{\Z}{{\mathbb Z}}

\newcommand{\Q}{{\mathbb Q}}
\newcommand{\R}{{\mathbb R}}
\newcommand{\K}{{\mathbb K}}
\newcommand{\pp}{\phi_{\K}}
\newcommand{\sk}{\mathcal{S}_k}
\newcommand{\Aa}{A_{G}(\vec{\rho},\vec{v},\vec{\lambda};d)}
\newcommand{\Ab}{A_{G}(\vec{\rho},\vec{v},\vec{\lambda};d)}
\newcommand{\Acb}{A_{G}(\vec{\rho},\vec{v},\vec{\lambda},\vec{\epsilon};d)}
\newcommand{\Aca}{A_{G}(\vec{\rho},\vec{v},\vec{\lambda},\vec{1};d)}
\newcommand{\A}{\mathcal{A}}

\newcommand{\DD}{\mathrm{Div}(k)}

\newcommand{\Ba}{B_{G}(\vec{\rho},\vec{v};d)}
\newcommand{\Bcb}{B_{G}(\vec{\rho},\vec{v},\vec{\epsilon};d)}
\newcommand{\Bca}{B_{G}(\vec{\rho},\vec{v},\vec{1};d)}
\newcommand{\Ccb}{F_{G}(\vec{\rho},\vec{v},\vec{\epsilon};d)}

\newcommand{\I}{{\mathbf I}}
\newcommand{\CC}{{\mathbb C}}

\newcommand{\MM}{m_{\rho,\K}}

\newcommand{\GG}{m_{\rho,\K}}

\newcommand{\F}{{\mathbb F}_q}

\newcommand{\FC}{\widetilde{E}_G(k,d)}

\newcommand{\ud}{\mathrm{d}}
\newcommand{\prob}{\mathrm{Prob}}
\newcommand{\fb}{\mathrm{Frob}}

\renewcommand{\vec}[1]{\underline{#1}}

\usepackage{setspace}
\begin{document}


\title[Zeta zeroes of abelian covers of algebraic curves]{Distribution of zeta zeroes for abelian covers of algebraic curves over a finite field}

\author[M. Xiong]{\sc Maosheng Xiong}


\thanks{The author was supported by the Research Grants Council of Hong Kong under Project Nos. RGC606211 and DAG11SC02.}


\begin{abstract}

For a function field $k$ over a finite field with $\F$ as the field of constant, and a finite abelian group $G$ whose exponent is divisible by $q-1$, we study the distribution of zeta zeroes for a random $G$-extension of $k$, ordered by the degree of conductors. We prove that when the degree goes to infinity, the number of zeta zeroes lying in a prescribed arc is uniformly distributed and the variance follows a Gaussian distribution.
\end{abstract}

\maketitle

\thispagestyle{empty}

\maketitle

\thispagestyle{empty}

\section{Introduction}

There has been some interest recently in statistics of zeroes of zeta functions for curves over a finite field $\F$. The fascinating theory of Katz and Sarnak (\cite{kat}) roughly states that as $q \to \infty$, the limiting distribution of zeroes of L-functions in a family of ``geometric objects'' defined over $\F$ is the same as that of eigenvalues of random matrices in certain monodromy groups. Here, however, the condition $q \to \infty$ is necessary as the argument depends on Deligne's equidistribution theorem (\cite{Del74}). One question hence remains: what happens if the field size $q$ is fixed and other parameters go to infinity?

In this direction Faifman and Rudnick (\cite{fai}) studied statistics of zeroes of zeta functions for the family of hyperelliptic curves of genus $g$ over $\F$. This is given explicitly by the affine model
\begin{eqnarray} \label{1:cf} C_f: y^2=f(x), \end{eqnarray}
where $f$ runs over the set of monic square-free polynomials of degree $2g+2$. The zeta function of $C_f/\F$ has $2g$ zeroes, and by the Riemann Hypothesis for curves (\cite{wei}), all of them lie on a circle of radius $q^{-1/2}$. They proved that as $g \to \infty$, the number of such zeta angles inside a fixed interval $\I \subset (-1/2,1/2)$ is asymptotically $2g|\I|$, here $|\I|$ is the length of $\I$, and the variance of this quantity as $C_f$ varies in the family is a Gaussian distribution; moreover, the result holds for shrinking intervals as long as $2g|\I|$ tends to infinity. This family form the moduli space of hyperelliptic curves of a fixed genus on which the distribution result can be reformulated. This result is consistent with the random matrix model predicted by the theory of Katz and Sarnak when both $g$ and $q$ tend to infinity.

This beautiful work of Faifman and Rudnick was extended to the family of $l$-fold covers of the projective line (\cite{xio}), and more recently to the family of Artin-Schreier covers of the projective line (\cite{buc3}), on which similar distribution results were obtained. The results can also be adapted to moduli spaces for these two families. Interested readers may refer to \cite{xio,buc3} for details. In the same spirit with respect to other statistics, in particular, the distribution of the number of $\F$-points on a family of curves or similar statistics as $g \to \infty$, this is slightly easier, and distribution results have been obtained for quite a few families (see for example \cite{buc1,buc11,buc2,bk,ent,kur,kur2,woo1}).


The above families of curves, as interesting and important as they are, could be interpreted in a different way. Let $k=\F(x)$ be the rational function field. Hyperelliptic curves $C_f$ defined in (\ref{1:cf}) correspond to quadratic extensions $k(\sqrt{f(x)})$ of $k$ with Galois group $\Z/2\Z$. Similarly $l$-fold covers of the projective line in \cite{xio} correspond to function field extensions of $k$ with Galois group isomorphic to $\Z/l\Z$ (here $q \equiv 1 \pmod{l}$), and Artin-Schreier curves considered in \cite{buc3} correspond to extensions of $k$ with Galois group isomorphic to $\Z/p\Z$, where $p$ is the characteristic of $\F$. Since the zeta function of a curve over $\F$ is the same as the zeta function of its function field (\cite{mor}), the results of \cite{fai,xio,buc3} can be summarized as distribution of zeroes of zeta functions for function fields running over abelian extensions of $k$ with the Galois group isomorphic to $\Z/2\Z,\Z/l\Z$ and $\Z/p\Z$ respectively. It is natural to investigate the question: how are the zeta zeros distributed as function fields run over abelian extensions of any function field with any fixed finite Galois group? Or in other words, how are the zeta zeros distributed for abelian covers of an algebraic curve over a finite field with a fixed Galois group?


From now on we fix an arbitrary function field $k$ with field of constant $\F$ and a finite abelian group $G$. Define a $G$-extension of $k$ to be a Galois extension $\K/k$ with an isomorphism $\pp:\mathrm{Gal}(\K/k) \to G$. An isomorphism of two $G$-extensions $\K$ and $\K'$ is given by an isomorphism $\K \to \K'$ of $k$-algebras that respects the $G$-action on $\K$ and $\K'$. Let $E_G(k)$ be the set of isomorphism classes of $G$-extensions of $k$. For a $G$-extension $\K$, let $\mathrm{Cond}(\K)$ be the conductor of $\K$ over $k$. For each positive integer $d$, we define
\begin{eqnarray} \label{1:ekg} E_G(k,d)=\{\K \in E_G(k): \deg_k \mathrm{Cond}(\K)=d\}. \end{eqnarray}
It will be clear that this is a finite set, with $\#E_G(k,d) \to \infty$ as $d \to \infty$. We assign the uniform probability measure on $E_G(k,d)$.

We shall study the distribution of zeroes for the zeta function $\zeta_{\K}(s)$ as $\K$ runs over the set $E_G(k,d)$ under the limit $d \to \infty$. It is known that $\left(1-q^s\right)\left(1-q^{1-s}\right) \zeta_{\K}(s)$ is a polynomial of degree $2g_{\K}$ in $q^{-s}$, where $g_{\K}$ is the genus of $\K$ and satisfies the Riemann hypothesis (\cite{ros}), so we may write
\[\zeta_{\K}(s)=\frac{\prod_{i=1}^{2g_{\K}} \left(1-\sqrt{q} e(\theta_{\K,i})q^{-s}\right)}{\left(1-q^s\right)\left(1-q^{1-s}\right)}, \]
here $e(\alpha)=e^{2 \pi i \alpha}$ and the angles satisfy $\{\theta_{\K,i}:i\} \subset [-\frac{1}{2},\frac{1}{2})$. Now we fix a subinterval $\I \subset (-\frac{1}{2},\frac{1}{2})$, and for simplicity as in \cite{fai}, we assume that $\I$ is symmetric around the origin. Define
\[N_{\I}(\K)=\#\left\{i: \theta_{\K,i} \in \I\right\}. \]
We will study how the quantity $N_{\I}(\K)$ is distributed as $\K$ runs over the set $E_G(k,d)$ with $d \to \infty$. We prove \begin{thm} \label{main2} Let $k$ be a function field with field of constant $\F$ and $G$ be a finite abelian group of order $\kappa$ such that $q \equiv 1 \pmod{\exp(G)}$ where $\exp(G)$ is the exponent of $G$. Let $\I \subset \left(-\frac{1}{2},\frac{1}{2}\right)$ be a symmetric subinterval and assume that $d |\I| \to \infty$ as $d \to \infty$ where $|\I|$ is the length of $\I$. Then for any real numbers $a,b$, we have
\begin{eqnarray*} \lim_{d \to \infty} \prob_{E_G(k,d)} \left(a<\frac{N_{\I}(\K)-g_{\K}|\I|}{\sqrt{\frac{2(\kappa-1)}{\pi^2} \log (g_{\K} |\I|)}}<b \right) =\frac{1}{\sqrt{2 \pi}} \int_{a}^{b} e^{-\frac{x^2}{2}} \, \ud x \,.\end{eqnarray*}
\end{thm}
We remark that the technical condition $q \equiv 1 \pmod{\exp (G)}$ effectively excludes Artin-Schreier extensions considered in \cite{buc3}. As was shown in \cite{buc3}, even over the rational function field, the treatment of Artin-Schreier extensions is quite different from \cite{fai,xio}, which means that to include such extensions in Theorem \ref{main2} might be technically complicated. We also remark that Theorem \ref{main2} is consistent with results in \cite{fai,xio,buc3} when $\# G=2,l \ge 3$ and $p$ respectively.

As was noted in \cite{xio} and \cite{buc3}, there is a subtle structure, that is,
\[\zeta_{\K}(s)=\zeta_k(s) \prod_{1 \ne \rho \in \widehat{G}} L(k,\rho \circ \pp,s),\]
where $L(k,\rho \circ \pp,s)$ is the L-function with respect to a non-principal character $\rho$. Since each $L(k,\rho \circ \pp,s)$ satisfies the Riemann hypothesis, it is more natural to study the distribution of zeroes of $L(k,\rho \circ \pp,s)$ for a fixed $\rho$ as $\K$ varies in $E_G(k,d)$. To be more precisely, for any non-principal character $\rho \in \widehat{G}$, it is known that each $L(k,\rho \circ \pp,s)$ is a polynomial of degree say $m_{\rho,\K}$ in $q^{-s}$ with constant term 1 and satisfies the Riemann hypothesis (\cite{ros}), so we may write
\[L(k,\rho \circ \pp,s)=\prod_{i=1}^{m_{\rho,\K}} \left(1-\sqrt{q} e(\theta_{\rho,\K,i})q^{-s}\right), \]
here $e(\alpha)=e^{2 \pi i \alpha}$ and the angles satisfy $\{\theta_{\rho,\K,i}:i\} \subset [-\frac{1}{2},\frac{1}{2})$. The collection of zeta angles $\left\{\theta_{\rho,\K,i}: 1 \le i \le m_{\rho,\K}\right\}$ for all the $\rho$'s (together with zeros of $\zeta_k(s)$) gives all the zeros of $\zeta_{\K}(s)$.  Now we fix a symmetric subinterval $\I \subset (-\frac{1}{2},\frac{1}{2})$, and define
\[N_{\rho,\I}(\K)=\#\left\{i: \theta_{\rho,\K,i} \in \I\right\}. \]
We will investigate a more subtle question, that is, what is the joint distribution of the quantities $N_{\rho,\I}(\K)$ for non-principal characters $\rho$ when $\K$ runs over the set $E_G(k,d)$ as $d \to \infty$. Since clearly $N_{\rho,\I}(\K)=N_{\rho^{-1},\I}(\K)$, from possibly two characters $\rho,\rho^{-1}$ we only need to consider one of them. We prove
\begin{thm} \label{main} Let $k$ be a function field with field of constant $\F$ and $G$ be a finite abelian group such that $q \equiv 1 \pmod{\exp(G)}$ where $\exp(G)$ is the exponent of $G$. Let $\rho_i: G \to \CC^*, 1 \le i \le t$ be distinct non-principal characters of $G$ such that $\rho_i \rho_j \ne 1$ for any $i \ne j$. Let $\I \subset \left(-\frac{1}{2},\frac{1}{2}\right)$ be a symmetric subinterval and assume that $d |\I| \to \infty$ as $d \to \infty$ where $|\I|$ is the length of $\I$. Then for any real numbers $a_i,b_i, 1 \le i \le l$, we have
\begin{eqnarray*} \lim_{d \to \infty} \prob_{E_G(k,d)} \left(a_i<\frac{N_{\rho_i,\I}(\K)-m_{\rho_i,\K}|\I|}{\sqrt{\frac{r_{\rho_i}}{\pi^2} \log (m_{\rho_i,\K} |\I|)}}<b_i, \,\, \forall \,\, 1 \le i \le t \right) =\prod_{i=1}^t \frac{1}{\sqrt{2 \pi}} \int_{a_i}^{b_i} e^{-\frac{x^2}{2}} \, \ud x \,,\end{eqnarray*}
where
\[r_{\rho}=\left\{\begin{array}{lll}
1&:& \mbox{ if the order of } \rho \mbox{ is } \ge 3;\\
2&:& \mbox{ if the order of } \rho \mbox{ is } 2.
\end{array}\right.\]
\end{thm}

\noindent {\bf Remark.} (a). Compared with Theorem \ref{main}, the number $r_{\rho}$ is missing in \cite[Theorem 1]{xio}. That was a typo. Actually the proof there was correct and clearly indicated that the factor $2$ should be replaced by $r_{\rho}$, which depends on whether or not $\rho^2=1$. This is consistent with \cite{fai,buc3}.

(b). The family of $G$-extensions of $k$ form the moduli space of abelian covers of an algebraic curve over $\F$ on which similar results as Theorems \ref{main2} and \ref{main} can formulated...

(c). We count $G$-extensions of $k$ by conductors (see the definition of $E_G(k,N)$ given in (\ref{1:ekg})). It might be more interesting to count $G$-extensions by discriminants, as discriminants determine the genus of the function fields (see the Riemann-Hurwitz formula (\cite[Theorem 7.16]{ros})). However, this turns out to be more difficult, as already noted in \cite{woo} (see also \cite{wri}), and we are not able to do it here. We might pursue it in the future.

The paper is organized as follows. In Section 2, we collect several results which will be used later. In Sections 3-4, we analyze certain sums over a larger set $\widetilde{E}_G(k,d)$ (see the definition in (\ref{0:egkd})). This part is inspired by the paper \cite{woo}. In Section 5 we prove Theorem \ref{main2} for the set $\widetilde{E}_G(k,d)$ as $d \to \infty$. In the final section, Section 6, Theorem \ref{main2} on $E_G(k,d)$ is derived from the results on $\widetilde{E}_G(k,d)$, and Theorem \ref{main} is also derived from Theorem \ref{main2}.


\section{Preliminaries}

In this section we collect several results which will be used later. Interested readers may refer to \cite{ros} for more details.

\subsection{L-functions and the explicit formulas}

Let $k$ be a function field with field of constant $\F$ and $\K/k$ be a Galois extension with an isomorphism $\pp: \mathrm{Gal}(\K/k) \to G$ where $G$ is a finite abelian group. Let $\rho: G \to \CC^*$ be a non-principal character. The L-function $L(k,\rho \circ \phi_K,s)$ for a complex variable $s$ is given by
\[L(k,\rho \circ \phi_K,s)=\prod_{v \in \sk} \left(1-\rho \circ \pp \left(v\right) |v|^{-s}\right)^{-1}, \]
where $\sk$ is the set of places of $k$, and for any $v \in \sk$, $|v|:=q^{\deg v}$, and $\rho \circ \pp(v):=\rho \circ \pp \left((v,\K/k)\right)$ where $(v,\K/k)$ is the Artin symbol when $v$ is unramified in $\rho \circ \pp$; if $v$ is ramified in $\rho \circ \pp$, we simply define $\rho \circ \pp\left(v\right):=0$. It is known from the Riemann hypothesis for curves that $L(k,\rho \circ \pp,s)$ is polynomial of finite degree $m_{\rho,\K}$ in $q^{-s}$ with constant term 1, so we can write
\[L(k,\rho \circ \pp,s)=\prod_{i=1}^{m_{\rho,\K}} \left(1-\sqrt{q}e(\theta_{\rho,\K,i})q^{-s}\right).\]
The degree $m_{\rho,\K}$ is given by the formula (\cite{ros}) \[m_{\rho,\K}=2g_k-2+\deg_k \mathcal{F}(\rho \circ \pp), \]
where $g_k$ is the genus of the field $k$ and $\mathcal{F}(\rho \circ \pp)$ is the Artin conductor of $\rho \circ \pp$. From this we see that
\[m_{\rho,\K} \asymp d \quad \mbox{ as } d=\deg_k \mathrm{Cond}(\K) \to \infty. \]
Taking logarithmic derivative of $L(k,\rho \circ \pp,s)$ with respect to $s$ by using the two different expressions and equating the coefficients, we obtain the the so-called ``explicit formulas''
\begin{eqnarray} \label{2:exp} \sum_{i=1}^{m_{\rho,\K}}e(n \theta_{\rho,\K,i})=-q^{-|n|/2}\sum_{\substack{v \in \sk \\
\deg v\mid n}} (\deg v) \rho \circ \pp(v)^{n/\deg v}, \,\, \forall \, 0 \ne n \in \Z. \end{eqnarray}

\subsection{Beurling-Selberg functions}

Let $\I=[-\beta/2,\beta/2]$ be an interval, symmetric about the origin, of length $0<\beta<1$, and $K \ge 1$ an integer. The Beurling-Selberg polynomials $I^{\pm}_K$ are trigonometric polynomials approximating the indicator function $\mathbf{1}_{\I}$ satisfying (see the exposition in \cite[Chapter 1.2]{mon}): \begin{itemize}
\item $I^{\pm}_K$ are trigonometric polynomials of degree $\le K$.

\item Monotonicity: $I^{-}_K \le \mathbf{1}_{\I} \le I^{+}_K$.

\item The integral of $I^{\pm}_K$ is close to the length of the interval:
    \begin{eqnarray} \label{2:IK1} \int_0^1 I^{\pm}_K(x) \, \ud x = \int_0^1 \mathbf{1}_{\I}(x) \, \ud x \pm \frac{1}{K+1}\end{eqnarray}

\item $I^{\pm}_K(x)$ are even (since the interval $\I$ is symmetric about the origin).
    \end{itemize}
As a consequence of (\ref{2:IK1}), the non-zero Fourier coefficients of $I^{\pm}_K(x)$ satisfy
\[\left|\widehat{I}^{\pm}_K(k) - \widehat{\mathbf{1}}_{\I}(k)\right| \le \frac{1}{K+1}\]
and in particular
\[\left|\widehat{I}^{\pm}_K(k)\right| \le \frac{1}{K+1}+ \min \left(\beta, \frac{\pi}{|k|}\right), \quad 0<|k| \le K. \]
This implies
\begin{eqnarray*} \left|\widehat{I}^{\pm}_K(k) k\right| \ll 1, \quad k \in \Z.
\end{eqnarray*}
\begin{itemize}
\item If $K\beta>1$, then (\cite[Propsition 4.1]{fai})
\begin{eqnarray*} \label{2:IK31} \sum_{n \ge 1}  \widehat{I}^{\pm}_K(2n)=O(1),\end{eqnarray*}
\begin{eqnarray} \label{2:IK3} \sum_{n \ge 1} n \widehat{I}^{\pm}_K(n)^2=\frac{1}{2 \pi^2} \log K\beta+O(1).\end{eqnarray}
\end{itemize}
All the implied constants above are independent of $K$ and $\beta$. We consider
\[\sum_{v \in \sk} \widehat{I}^{\pm}_K(\deg v)^2 (\deg v)^2 |v|^{-1}.\]
The prime number theorem $\# \left\{v \in \sk:\deg v=n\right\} =q^n/n+O\left(q^{n/2}\right)$ gives
\begin{eqnarray*} \sum_{v \in \sk} \widehat{I}^{\pm}_K(\deg v)^2 (\deg v)^2 |v|^{-1}&=&\sum_{1 \le n \le K} \widehat{I}^{\pm}_K(n)^2 \left(n+O\left(nq^{-n/2}\right)\right). \end{eqnarray*}
Using (\ref{2:IK3}) we obtain (similar to Equation (7.4) in \cite{fai})
\[\sum_{v \in \sk} \widehat{I}^{\pm}_K(\deg v)^2 (\deg v)^2 |v|^{-1}=\frac{1}{2 \pi^2} \log K\beta +O(1)\,. \]

\section{Abelian extensions with a fixed Galois group: part I}

\noindent {\bf Notation.} Let $k$ be a function field with field of constant $\F$ and $G$ be a finite abelian group of order $\kappa$, given explicitly by $G=\prod_{i=1}^t \Z/n_i\Z, \, n_{i+1}\mid n_i$, hence $\exp (G)=n_1$. We assume that $q \equiv 1 \pmod{n_1}$. For all positive integers $n$ such that $\gcd(n,q)=1$, we choose compatible systems of primitive $n$-th roots of unity $\{\xi_n\} \subset \bar{k} \hookrightarrow \bar{k}_v$ and $\{\zeta_n\} \subset \bar{\Q}$ such that if $n' \mid n$, then $\xi_{n'}=\xi_n^{n/n'}$ and $\zeta_{n'}=\zeta_n^{n/n'}$. Let $\iota$ be the map such that $\iota(\xi_n)=\zeta_n, \, \gcd(n,q)=1$.
\[\xymatrix{
&{\bar{k}} \ar@{-}[d]\ar@{-}[r] & {{\bar{k}_v}} \ar@{-}[d] & {\{\xi_n\}} \ar[r]^{\iota} &{\{\zeta_n\}} & {\bar{\Q}} \ar@{-}[d] \\
{\F} \ar@{-}[r] &{k} \ar@{-}[r]& {k_v}
& {}                                             & {} & {\Q}}\]
Define
\[\widetilde{E}_G(k,d):= \bigcup_{H \,\, \mbox{\tiny subgroup of } \, G}E_H(k,d), \]
or more precisely
\begin{eqnarray} \label{0:egkd} \widetilde{E}_G(k,d):= \left\{\K: \mathrm{Gal}(\K/k) \overset{\pp}{\longrightarrow} G \, \mbox{ injective}, \deg_k \mathrm{Cond}(\K)=d\right\}. \end{eqnarray}
The purpose of this section is to prove asymptotic formulas on $\widetilde{E}_G(k,d)$ which will be used later.

\begin{thm} \label{3:rho} For any characters $\rho_1,\ldots,\rho_r \in \widehat{G}$ of order say $\tau_1,\ldots,\tau_r$ respectively, any distinct places $v_1,\ldots,v_r \in \sk$ and any integers $\lambda_1,\ldots,\lambda_r$, define $\vec{\rho}=(\rho_1,\ldots,\rho_r),\vec{v}=(v_1,\ldots,v_r),
\lambda=(\lambda_1,\ldots,\lambda_r)$ and consider the sum
\begin{eqnarray} \label{3:A}
\Aa:=\sum_{\K \in \widetilde{E}_G(k,d) } \prod_{i=1}^r \rho_i \circ \pp\left(v_i\right)^{\lambda_i}.
\end{eqnarray}
There is a constant $C=C(G,k)$ depending only on $G$ and $k$ such that

\noindent {\bf (i).} if $\tau_i \nmid \lambda_i$ for some $i$, then
\[ \Aa \ll_{G,k,\eta} C^{r} q^{(1+\eta)d/{2}};\]

\noindent {\bf (ii).} if $\tau_i \mid \lambda_i$ for any $1 \le i \le r$, then
\[\Aa=H \cdot  H_{_{\sum_0}}c_k^{\kappa-1}d^{\kappa-2} q^{d+\kappa-1} \left\{1+O_{G,k}\left(C^rd^{-1}\right)\right\}. \]
Here $H$ is an absolute constant given by
\[H=\prod_{v \in \sk}\left(1+(\kappa-1)|v|^{-1}\right)
\left(1-|v|^{-1}\right)^{\kappa-1},\]
and the constant $H_{_{\sum_0}}$ is given by
\[H_{_{\sum_0}}=\prod_{i=1}^r \frac{1+(\#\ker \rho_i-1)|v_i|^{-1}}{
1+(\kappa-1)|v_i|^{-1}}.\]
The constant $c_k$ depends only on $k$ and can be given explicitly
\[c_k={(q-1)^{-1}q^{-g_k}}h_k,\]
where $h_k$ is the class number of $k$ and $g_k$ is the genus of $k$.

\end{thm}

To prove Theorem \ref{3:rho}, we use class field theory together with techniques borrowed from \cite{woo,wri}. Actually for (ii) of Theorem \ref{3:rho}, the quantity $\Aa$ counts the number of abelian extensions $\K \in \widetilde{E}_G(k,d)$ such that $v_i$ is unramified in $\rho_i \circ \pp$ for each $i$. In particular if $r=0$, then (ii) of Theorem \ref{3:rho} implies that
\begin{eqnarray} \label{3:egkd} \#\widetilde{E}_G(k,d)=H \cdot \binom{d+\kappa-2}{\kappa-2}  c_k^{\kappa-1}q^{d+\kappa-1} \left\{1+O_{G,k}\left(C^rd^{-1}\right)\right\}. \end{eqnarray}
This is the function field analogue of results in \cite{woo,wri} with explicit error terms. It is conceivable that the class field theory and techniques from \cite{woo,wri} can be adapted to the function field setting to obtain asymptotic formulas on $\Aa$ without much difficulty. However, the main issue here is to obtain error terms strong enough for applications in mind. As is already clear in \cite{woo,wri}, and as we shall see, this is not easy. We are actually quite lucky to get a power saving as in (i) of Theorem \ref{3:rho}. If, for example, the function fields were ordered by discriminants (i.e., the set $\widetilde{E}_G(k,d)$ consists of such $\K$'s with $\deg_k \mathrm{disc}(\K)=d$), asymptotic formulas for $\Aa$ can still be obtained, however, the error terms were too weak to be of any use for the purpose of applications. We might try to tackle this problem in the future.

Now we start the proof of Theorem \ref{3:rho}. We follow colsely the arguments in \cite{woo}. For information on class field theory, interested readers may refer to the paper \cite{wri} for a practical overview and the book \cite{neu} for a comprehensive treatment of the theory.

\subsection{Preparation}

By class field theory, abelian extensions $\K \in \widetilde{E}_G(k,d)$ corresponds one-to-one to homomorphisms $\chi: J/k^* \to G$ with $\mathrm{Cond}(\chi)=\mathrm{Cond}(\K)$ where $J$ is the group of id\`{e}les of $k$ and $J/k^*$ is the id\`{e}le class group. Here $k^*$ is understood to be its embedding image in $J$. For simplicity let us define
\[C(\chi):=| \mathrm{Cond}(\chi) |=q^{\deg_k \mathrm{Cond}(\K)}. \]
For each place $v$ of $k$, let $k_v,o_v$ be the local field and the local ring at $v$ respectively. Denote by $o_v^*$ the group of units in $o_v$. We shall identify $k_v^*$ and $o_v^*$ with their images in $J$. Let $\pi_v \in o_v$ be a generator of the unique maximal ideal of $o_v$, then $\pi_v$ can also be regarded as an element in $k_v^* \subset J$. For $\chi: J/k^* \to G$, let $\chi_v:k_v^*k^*/k^* \to G$ be the $v$-th component of $\chi$. Since $q \equiv 1 \pmod{n_1}$, if $v$ is ramified then $(1+\pi_v o_v)k^*/k^* \subset \ker (\chi_v)$. Define
\[c(\chi_v):=\left\{\begin{array}{ccl}\
1&:& v \mbox{ is ramified}; \\
0&:& v \mbox{ is unramified}.
\end{array}\right.\]
Then we have
\[C(\chi)=\prod_{v \in \sk}N v^{c(\chi_v)}, \]
where $Nv:=|v|=q^{\deg_k v}$ is the absolute norm of $v$. With this preparation and by class field theory we can rewrite $\Aa$ in (\ref{3:A}) as \[\Aa=\sum_{\substack{\chi:J/k^* \to G \\
v_i \, \mbox{\tiny unramified in } \rho_i \circ \chi \, \forall \, i\\
C(\chi)=q^d
}} \prod_{i=1}^r \rho_i \circ \chi (\pi_{v_i})^{\lambda_i}. \]
Let
\[ \Sigma_0:=\{v_1,\ldots,v_r\},\]
and we choose a finite subset $\Sigma' \subset \sk$ of order say $c$ such that $o_{\Sigma'}$, the ring of $\Sigma'$-integers of $k$, has class number 1. Define
\[\Sigma=\Sigma_0 \bigcup \Sigma', \quad r'=\#\Sigma \le r+c. \]
Then the class number of $o_{\Sigma}$ is also 1, and the natural map $J_{\Sigma}/o_{\Sigma}^* \to J/k^*$ is an isomorphism (see \cite[Lemma 2.8]{woo}), here $J_{\Sigma}$ is the group of id\`{e}les which have components in $o_v^*$ for all places $v \notin \Sigma$, and $o_{\Sigma}^*$ is the group of units in $o_{\Sigma}$. Hence we have
\begin{eqnarray*}
\Aa=\sum_{\substack{\chi:J_{\Sigma}/o_{\Sigma}^* \to G \\
v_i \, \mbox{\tiny unramified in } \rho_i \circ \chi \, \forall \, i\\
C(\chi)=q^d}} \prod_{i=1}^r \rho_i \circ \chi (\pi_{v_i})^{\lambda_i}.
\end{eqnarray*}
Let $\A=\prod_{i=1}^t o_{\Sigma}^*/o_{\Sigma}^{*n_i}$. Given a $\chi: J_{\Sigma} \to G$ with projection $\chi_i: J_{\Sigma} \to \Z/n_i\Z$ (or the same from $k_v^*$ or $o_v^*$), and an $\vec{\epsilon}=(\epsilon_1,\ldots,\epsilon_k) \in \A$, we define $\dot{\chi}(\vec{\epsilon})=\prod_{i=1}^t
\zeta_{n_i}^{\chi_i(\epsilon_i)}$, where we evaluate $\chi_i(\epsilon_i)$ using the natural map $o_{\Sigma}^* \to J_{\Sigma}$ (or to $k_v^*$ or $o_{v}^*$). We define the twists
\begin{eqnarray} \label{3:twi} \Acb=
\sum_{\substack{\chi:J_{\Sigma} \to G \\
v_i \, \mbox{\tiny unramified in } \rho_i \circ \chi \, \forall \, i\\
C(\chi)=q^d}} \left\{\prod_{i=1}^r \rho_i \circ \chi (\pi_{v_i})^{\lambda_i}\right\} \cdot \dot{\chi}(\vec{\epsilon}). \end{eqnarray}
It is known from \cite[Corollary 2.9]{woo} that
\begin{eqnarray} \label{3:epf} \Ab=\frac{1}{\#\A}
\sum_{\vec{\epsilon} \in \A}\Acb. \end{eqnarray}
Thus to find the asymptotic behavior of $\Ab$ as $d \to \infty$, it suffices to study $\Acb$ for each $\vec{\epsilon}$. A standard technique is to analyze the generating function
\begin{eqnarray*} \label{3:gb}
F_{G}(\vec{\rho},\vec{v},\vec{\lambda},\vec{\epsilon};s):=\sum_{\substack{\chi:J_{\Sigma} \to G \\
v_i \, \mbox{\tiny unramified in } \rho_i \circ \chi \, \forall \, i}} \frac{\prod_{i=1}^r \rho_i \circ \chi (\pi_{v_i})^{\lambda_i}}{C(\chi)^s} \cdot \dot{\chi}(\vec{\epsilon}),
\end{eqnarray*}
here $s$ is a complex variable. The functions $F_G(\vec{\rho},\vec{v},\vec{\lambda},\vec{\epsilon};s)$ are convenient to work with because they have Euler products (\cite{woo,wri})
\begin{eqnarray*} F_G(\vec{\rho},\vec{v},\vec{\lambda},\vec{\epsilon};s)&=&
\prod_{v \notin \Sigma }\left(\sum_{\chi_v:o_v^* \to G} \frac{\dot{\chi}_v(\vec{\epsilon})}{Nv^{c(\chi_v)s}}\right) \prod_{v \in \Sigma \setminus \Sigma_0}\left(\sum_{\chi_v:k_v^* \to G} \frac{\dot{\chi}_v(\vec{\epsilon})}{Nv^{c(\chi_v)s}}\right) \times \\
&& \prod_{i=1 }^r\left(\sum_{\substack{\chi_{v_i}:k_{v_i}^* \to G\\
v_i \, \mbox{\tiny unramified in } \rho_i \circ \chi_{v_i}}} \frac{\rho_i \circ \chi_{v_i}(\pi_{v_i})^{\lambda_i} \dot{\chi}_{v_i}(\vec{\epsilon})}{Nv_i^{c(\chi_{v_i})s}}
\right).\end{eqnarray*}


\subsection{Analysis of the Euler products}

To study analytic behavior of $F_G(\vec{\rho},\vec{v},\vec{\lambda},\vec{\epsilon};s)$ for the complex variable $s$, we need to analyze the Euler factor at each place $v$, in particular the values $\dot{\chi}_v(\vec{\epsilon})$. For the number field case, this was done by \cite[Lemmas 2.15 and 2.16]{woo}, which were motivated by the work of Taylor (\cite{tay}). The function field analogues of the two lemmas follow literally the same lines of argument, so we only record the results here.

For each place $v$, we choose a generator $y_v$ of the tame inertia group of $k_v$ (which is isomorphic to $(o_v/v)^*$) such that $y_v \equiv \xi_{Nv-1} \pmod{v}$, where $\xi_{Nv-1}$ is the primitive $(Nv-1)$-th root of unity in $\bar{k} \subset \bar{k}_v$ which we have fixed at the beginning of this section. Then similar to \cite[Lemma 2.15]{woo} we have
\[\xi_{Nv-1}=\frac{\fb_v\left(y_v^{1/(Nv-1)}
\right)}{y_v^{1/(Nv-1)}}, \]
where the Frobenius is in the Galois group of the maximal unramified extension of $k_v$. \cite[Lemma 2.16]{woo} also applies in the function field setting, using our choice of the compatible systems of roots of unity $\{\xi_n\}$ and $\{\zeta_n\}$. The results become much easier because $q \equiv 1 \pmod{n_1}$, hence $\xi_{n_1} \in k$. We record the values $\dot{\chi}_v(\vec{\epsilon})$ as the following.
\begin{lemma} \label{3:chi}
For each place $v \notin \Sigma_0$, let $\chi_v(y_v)=g \in G$. Suppose the projections of $g$ to $\Z/n_i\Z$ are $n_ik_i/l_i \in \Z/n_i\Z$ where $l_i \mid n_i$ and $\gcd(k_i,l_i)=1$. Let $\vec{\epsilon}^g$ be the notation for $\prod_{i=1}^t \epsilon_i^{k_i/l_i}$. Then
\[\dot{\chi}_v(\vec{\epsilon})=\iota\left(\prod_{i=1}^{t}
\frac{\fb_v\left(\epsilon_i^{k_i/l_i}\right)}{\epsilon_i^{k_i/l_i}}\right)
=\iota\left(\frac{\fb_v\left(\vec{\epsilon}^g\right)}{\vec{\epsilon}^g}\right),\]
where the Frobenius is in the Galois group of the maximal extension of $k$ unramified outside $\Sigma_0$.
\end{lemma}
We observe
\begin{lemma} \label{3:chi2} For each $v \notin \Sigma_0$,
\[\sum_{\chi_v: o_v^* \to G} \dot{\chi}_v(\vec{\epsilon})=
 \left\{\begin{array}{ccc}
 \kappa&:& \vec{\epsilon}=\vec{1}; \\
 0&:& \vec{\epsilon} \ne \vec{1}. \end{array}\right.\]
\end{lemma}
\noindent {\bf Proof.} Let $A$ be the sum on the left. Since $q \equiv 1 \pmod{n_1}$ and $o_v^* \simeq (o_v/v)^*  \times \left(1+\pi_v o_v\right)$, the homomorphism $\chi_v: o_v^* \to G$ factors through $(o_v/v)^*$ and hence is completely determined by its value at the generator $y_v$, in other words, $\chi_v$ is determined by $g \in G$ such that $\chi_v(y_v)=g$. From Lemma \ref{3:chi} we have
\[A=\sum_{g \in G} \iota\left(\frac{\fb_v\left(\vec{\epsilon}^g\right)}{\vec{\epsilon}^g}\right). \]
If $\vec{\epsilon}=\vec{1}$, then clearly $A=\#G=\kappa$. Suppose $\vec{\epsilon} \ne \vec{1}$, then from the definition, there is a $g_0 \in G$ such that $\vec{\epsilon}^{g_0} \notin k$, hence $\frac{\fb_v\left(\vec{\epsilon}^{g_0}\right)}{\vec{\epsilon}^{g_0}} \ne 1$. So we have
\[\iota\left(\frac{\fb_v\left(\vec{\epsilon}^{g_0}\right)}{\vec{\epsilon}^{g_0}}\right) \cdot A=
\sum_{g \in G}\iota\left(\frac{\fb_v\left(\vec{\epsilon}^{g+g_0}\right)}{\vec{\epsilon}^{g+g_0}}\right)=A,\]
from which we derive that $A=0$. \quad $\square$

From Lemma \ref{3:chi2} we find that for any $ v \notin \Sigma$,
\begin{eqnarray} \label{3:chi3} \sum_{\chi_v: o_v^* \to G} \frac{\dot{\chi}_v(\vec{\epsilon})}{Nv^{c(\chi_v)s}}=
 \left\{\begin{array}{lcc}
1+(\kappa-1)|v|^{-s}&:& \vec{\epsilon}=\vec{1}; \\
1-|v|^{-s}&:& \vec{\epsilon} \ne \vec{1}. \end{array}\right.\end{eqnarray}
As for the Euler factor at $v \in \Sigma \setminus \Sigma_0$, using the isomorphism $k_v^* \simeq \langle \pi_v \rangle \times (o_v/v)^* \times \left(1+\pi_v o_v\right)$ and $q \equiv 1 \pmod{n_1}$, a homomorphism $\chi_v:k_v^* \to G$ is determined by the values $\chi_v(\pi_v)$ and $\chi_v(y_v)$. So
\[\sum_{\chi_v: k_v^* \to G} \frac{\dot{\chi}_v(\vec{\epsilon})}{Nv^{c(\chi_v)s}}=\sum_{\mu_v: \langle \pi_v \rangle \to G} \dot{\mu}_v(\vec{\epsilon})
\sum_{\gamma_v: o_v^* \to G} \frac{\dot{\gamma}_v(\vec{\epsilon})}{Nv^{c(\gamma_v)s}}. \]
The second sum on the right over $\gamma_v$ can be evaluated by (\ref{3:chi3}). As for the sum over $\mu_v$, by definition we find that
\[\sum_{\mu_v: \langle \pi_v \rangle \to G} \dot{\mu}_v(\vec{\epsilon}) = \prod_{i=1}^t \sum_{\mu_{v,i}: \langle \pi_v \rangle \to \Z/n_i\Z} \zeta_{n_i}^{\mathrm{ord}_v(\epsilon_i) \mu_{v,i}(\pi_v)}=\kappa \cdot \mathbf{1}_{n_i \mid \mathrm{ord}_v(\epsilon_i) \, \forall i}, \]
here $\mathbf{1}_{n_i \mid \mathrm{ord}_v(\epsilon_i) \, \forall i}=1$ if indeed $n_i \mid \mathrm{ord}_v(\epsilon_i) \, \forall 1 \le i \le t$; otherwise, the value is zero.

Finally for the Euler factor at $v_i \in \Sigma_0$, $1 \le i \le r$, which we denote by $A_{v_i}(s)$, the place $v_i$ is unramified in $\rho_i \circ \chi_{v_i}$ if and only if $o_{v_i}^* \subset \ker \rho_i \circ \chi_{v_i}$. Writing $\chi_{v_i}=(\mu_i, \gamma_i)$ where $\mu_i: \langle \pi_{v_i} \rangle \to G$ and $\gamma_i: (o_v/v)^* \to G$, we find
\begin{eqnarray} \label{3:chi4} A_{v_i}(s)=\sum_{\mu_{i}:\langle \pi_{v_i} \rangle \to G} \rho_i \circ \mu_i (\pi_{v_i})^{\lambda_i} \dot{\mu}_i(\vec{\epsilon})
\sum_{\substack{\gamma_i: (o_{v_i}/v_i)^* \to \ker \rho_i}} \frac{ \dot{\gamma}_{i}(\vec{\epsilon})}{Nv_i^{c(\gamma_i)s}}. \end{eqnarray}
If $\vec{\epsilon}=\vec{1}$, then
\[\sum_{\mu_{i}:\langle \pi_{v_i} \rangle \to G} \rho_i \circ \mu_i (\pi_{v_i})^{\lambda_i} =\left\{\begin{array}{ccc}
\kappa&:& \tau_i \mid \lambda_i,\\
0&:& \tau_i \nmid \lambda_i, \end{array}\right.\]
where $\tau_i$ is the order of $\rho_i$, and the second sum on the right of (\ref{3:chi4}) is
\[\sum_{\substack{\gamma_i: (o_{v_i}/v_i)^* \to \ker \rho_i}} \frac{ 1}{Nv_i^{c(\gamma_i)s}}=1+\frac{\# \ker \rho_i-1}{|v_i|^{s}}. \]
If $\vec{\epsilon} \ne \vec{1}$, we are contented with the fact that
\[|A_{v_i}(s)| \le \kappa \left(1+\frac{\# \ker \rho_i-1}{|v_i|^{\Re s}}\right).\]
Now we summarize the above analysis on $F_G(\vec{\rho},\vec{v},\vec{\lambda},\vec{\epsilon};s)$ as follows.

\begin{lemma} \label{3:fg} If $\vec{\epsilon}=\vec{1}$, and
\begin{itemize}
\item[(1.1).] if $\tau_i \mid \lambda_i$ for any $1 \le i \le r$, then
\[F_G(\vec{\rho},\vec{v},\vec{\lambda},\vec{1};s)=\kappa^{r'} \prod_{v \notin \Sigma_0} \left(1+\frac{\kappa-1}{|v|^{s}}\right) \prod_{i=1}^r\left(1+\frac{\# \ker \rho_i-1}{|v_i|^{s}}\right); \]

\item[(1.2).] if $\tau_i \nmid \lambda_i$ for some $i$, then \[F_G(\vec{\rho},\vec{v},\vec{\lambda},\vec{1};s)=0. \]
\end{itemize}
If $\vec{\epsilon} \ne \vec{1}$, then
\[F_G(\vec{\rho},\vec{v},\vec{\lambda},\vec{\epsilon};s)= h(s) \prod_{v \notin \Sigma_0} \left(1-\frac{1}{|v|^{s}}\right) ,  \]
where
\[|h(s)| \le \kappa^{r'} \prod_{i=1}^r \left(1+\frac{\# \ker \rho_i-1}{|v_i|^{\Re s}}\right). \]
\end{lemma}

\subsection{Two lemmas}
We now use Lemma \ref{3:fg} along with the function field version of the Tauberian Theorem (\cite[Chap 17]{ros}) to prove two estimates on $\Acb$ defined in (\ref{3:twi}).

\begin{lemma} \label{3:ad1}
If $\vec{\epsilon} \ne \vec{1}$, then
\[\Acb \ll_{G,k,\eta} C^r q^{(1+\eta)d/2}, \]
where $C=C(G,k) >0$ is a constant depending only on $G$ and $k$.
\end{lemma}

\noindent {\bf Proof.} For $\vec{\epsilon} \ne \vec{1}$,  from Lemma \ref{3:fg} we can write
\begin{eqnarray} \label{3:fga} F_G(\vec{\rho},\vec{v},\vec{\lambda},\vec{\epsilon};s)=\zeta_k(s)^{-1}H(s), \end{eqnarray}
where $\zeta_k(s)$ is the zeta function of $k$ defined by
\[\zeta_k(s)=\prod_{v \in \sk}\left(1-\frac{1}{|v|^s}\right)^{-1},\]
and the function $H(s)$ satisfies
\[|H(s)| \le \kappa^{r'} \prod_{i=1}^r \left(1+O(\kappa |v_i|^{-\Re s})\right). \]
Here $r'=\#\Sigma \le r+c$ where $c=\#\Sigma'$ which depends only on $k$. This shows that $F_G(\vec{\rho},\vec{v},\vec{\lambda},\vec{\epsilon};s)$ is analytic for all $s \in \CC$ except possible poles at $\Re s=\frac{1}{2}$. This analytic statement is equivalent to Lemma \ref{3:ad1} by using the standard Tauberian argument. For the function field version, it is most convenient to work on the variable $T=q^{-s}$. With no ambiguity we shall write the functions in $T$ as $\zeta_k(T), H(T)$ and $F_G(\vec{\rho},\vec{v},\vec{\lambda},\vec{\epsilon};T)$. It is known that $\zeta_k(T)$ is of the form
\begin{eqnarray} \label{3:zeta} \zeta_k(T)=\frac{\prod_{i=1}^{2g}\left(1-\sqrt{q} e(\theta_i) T\right)}{(1-T)(1-qT)}, \end{eqnarray}
where $g$ is the genus of $k$ and $\theta_i$'s are some real numbers. As for $H(T)$ we have
\begin{eqnarray} \label{3:ht} |H(T)| \le \kappa^{r'} \prod_{i=1}^r \left(1+O(\kappa |T|^{\deg v_i})\right). \end{eqnarray}
We may expand
\[F_G(\vec{\rho},\vec{v},\vec{\lambda},\vec{\epsilon};T)=
\sum_{d=0}^{\infty}T^d \Acb.\]
Since  $F_G(\vec{\rho},\vec{v},\vec{\lambda},\vec{\epsilon};T)$ is analytic for $T$ in the region $0<|T|<\tau=q^{-(1+\eta)/2}$, we find
\[\Acb=\frac{1}{2 \pi i} \oint_{|T|=\tau} \frac{F_G(\vec{\rho},\vec{v},\vec{\lambda},\vec{\epsilon};T)}{T^{d+1}} \, \ud  T,\]
where $\eta>0$ is an arbitrarily small real number. Now Lemma \ref{3:ad1} is immediate by estimating this integral, using (\ref{3:fga}) together with (\ref{3:zeta}) and (\ref{3:ht}). \quad $\square$

\begin{lemma} \label{3:ad2} If $\tau_i \mid \lambda_i$ for each $1 \le i \le r$, then
\[\Aca=\kappa^{r'} H \cdot H_{_{\sum_0}}c_k^{\kappa-1}d^{\kappa-2} q^{d+\kappa-1} \left\{1+O_{G,k}\left(C^rd^{-1}\right)\right\}, \]
where the constants $H,H_{\Sigma_0},c_k,C$ are the same as in Theorem \ref{3:rho}.
\end{lemma}

\noindent {\bf Proof.} From Lemma \ref{3:fg} we have
\begin{eqnarray*} \label{3:fs1} f(s):=F_G(\vec{\rho},\vec{v},\vec{\lambda},\vec{1};s)=\kappa^{r'} \zeta_k(s)^{\kappa-1} \widetilde{H}(s),\end{eqnarray*}
where
\begin{eqnarray} \label{3:hs2} \widetilde{H}(s)=\left\{ \prod_{v \in \sk} \left(1+\frac{\kappa-1}{|v|^s}\right)
\left(1-\frac{1}{|v|^s}\right)^{\kappa-1} \right\}  \left\{\prod_{i=1}^r\frac{\left(1+\frac{\# \ker \rho_i-1}{|v_i|^s}\right)}{\left(1+
\frac{\kappa-1}{|v_i|^s}\right)}\right\}. \end{eqnarray}
It is convenient to use $T=q^{-s}$ and we write the new functions on $T$ as $f,\zeta_k, \widetilde{H}$ as well. Now $f(T)$ has a pole at $T=q^{-1}$ of order $\kappa-1$, coming from $\zeta_k(T)^{\kappa-1}$, and $\widetilde{H}(T)$ is analytic for $T$ in the region $0<|T|<\tau=q^{-(1+\eta)/2}$ for an arbitrarily small real number $0<\eta<1$. We have
\[\frac{1}{2 \pi i} \oint_{|T|=\tau} \frac{f(T)}{T^{d+1}}\, \ud T=\mathrm{Res}_{T=q^{-1}}\left(\frac{f(T)}{T^{d+1}}\right)+
\mathrm{Res}_{T=0}\left(\frac{f(T)}{T^{d+1}}\right). \]
From the expression for $\zeta_k(T)$ in (\ref{3:zeta}) and $\widetilde{H}(s)$ in (\ref{3:hs2}), we find easily that
\[\frac{1}{2 \pi i} \oint_{|T|=\tau} \frac{f(T)}{T^{d+1}}\, \ud T \ll_{G,k} C^{r} q^{(1+\eta)d/2}. \]
On the other hand, expanding $f(T)$ as power series in $T$ we have
\[\Aca=\mathrm{Res}_{T=0}\left(\frac{f(T)}{T^{d+1}}\right). \]
Hence to find the asymptotics for $\Aca$, it suffices to compute the residue of $f(T)/T^{d+1}$ at $T=q^{-1}$, which we denote by $B$. The computation is straightforward: we use the formula
\[B=\frac{1}{(\kappa-2)!} \lim_{T \to q^{-1}} \left( (T-q^{-1})^{\kappa-1} \frac{f(T)}{T^{d+1}} \right)^{(\kappa-2)},\]
where the exponent $(\kappa-2)$ means to take the $(\kappa-2)$-th derivative with respect to $T$. Since $f(T)=\kappa^{r'} \zeta_k(T)^{\kappa-1} \widetilde{H}(T)$, expanding the derivative we obtain
\begin{eqnarray*} B&=&\frac{\kappa^{r'}}{(\kappa-2)!} \sum_{\substack{u,v,w \ge 0\\
u+v+w=\kappa-2}} \binom{\kappa-2}{u,v,w}(-1)^u(d+1) \cdots (d+u) q^{d+1+u} \times \\
&&\lim_{T \to q^{-1}}\widetilde{H}(T)^{(v)} \lim_{T \to q^{-1}} \left\{ (T-q^{-1})^{\kappa-1} \zeta_k(T)^{\kappa-1} \right\}^{(w)}. \end{eqnarray*}
The term from $u=\kappa-2,v=w=0$ gives the main contribution, which is
\begin{eqnarray*} B_0 &=& \frac{\kappa^{r'}}{(\kappa-2)!} (-1)^{\kappa-2}(d+1) \cdots (d+\kappa-2) q^{d+\kappa-1} \widetilde{H}(q^{-1}) (-c_k)^{\kappa-1}\\
&=& -\kappa^{r'} H \cdot H_{_{\sum_0}}c_k^{\kappa-1} d^{\kappa-2} q^{d+\kappa-1}\left(1+O_{G,k}(d^{-1})\right), \end{eqnarray*}
where the constants $H$ and $H_{\Sigma_0}$ come from (\ref{3:hs2}) with $s=1$ (i.e. $T=q^{-1}$), and
\[c_k=-\lim_{T \to q^{-1}} (T-q^{-1}) \zeta_k(T). \]
Using the expression (\cite[Theorem 5.9, page 53]{ros})
\[\zeta_k(T)=\sum_{n=0}^{2g-2} b_nT^n+ \frac{h_k}{q-1}\left(\frac{q^{g_k}}{1-qT}-\frac{1}{1-T}\right)
T^{2g-1},\]
where $h_k$ is the class number of $k$ and $g_k$ be the genus, we find easily
\[c_k=(q-1)^{-1}q^{-g_k}\, h_k.\]
All other terms in $B$ from $u <\kappa-1$ are bounded by
\[\ll_{G,k} C^r d^{-1}B_0.  \]
Combining the above estimates completes the proof of Lemma \ref{3:ad2}. \quad $\square$

\subsection{Proof of (i) and (ii) in Theorem \ref{3:rho}}
We are now ready to prove (i) and (ii) in Theorem \ref{3:rho}. If $\tau_i \nmid \lambda_i$ for some $i$, then for any $\vec{\epsilon} \ne \vec{1}$, from Lemma \ref{3:ad1} we have
\begin{eqnarray} \label{4:i} \Acb \ll_{G,k,\eta} C^{r} q^{(1+\eta)d/2}. \end{eqnarray}
From (1.2) in Lemma \ref{3:fg} we also know that
\[\Aca =0. \]
Then from (\ref{3:epf}) we conclude that
\[\Ab \ll_{G,k,\eta} C^{r} q^{(1+\eta)d/2}. \]
So (i) is proved. Now suppose $\tau_i \mid \lambda_i$ for any $1 \le i \le r$. For $\vec{\epsilon} \ne \vec{1}$, we still have the estimate (\ref{4:i}) as before. Hence from (\ref{3:epf}) we have
\[\Ab=\frac{\Aca}{\#\A} +O_{G,k,\eta}\left(C^{r} q^{(1+\eta)d/2}\right).
\]
The $\Aca$ can be read off from Lemma \ref{3:ad2}. Hence to prove (ii) of Theorem \ref{3:rho}, we need to show that $\#\A=\kappa^{r'}$. This is indeed the case: since $o_{\Sigma}^*/\F^*$ is a free group on $r'-1$ generators (\cite[Proposition 14.2]{ros}) and $q \equiv 1 \pmod{n_i}, \forall i$, we have
\[o_{\Sigma}^*/o_{\Sigma}^{*n_i} \simeq \F^*/\F^{*n_i} \times \left(\Z/n_i\Z\right)^{r'-1} \simeq \left(\Z/n_i\Z\right)^{r'}. \]
Therefore
\[\#\A=\prod_{i=1}^t \#\left(\Z/n_i\Z\right)^{r'}=\prod_{i=1}^t n_i^{r'}=\kappa^{r'}. \]
Now (ii) is proved. This completes the proof of Theorem \ref{3:rho}. \quad $\square$


\section{Abelian extensions with a fixed Galois group: part II}

Using notation from the previous section, we prove
\begin{thm} \label{4:rho} For any non-principal characters $\rho_1,\ldots,\rho_r \in \widehat{G}$ and any distinct places $v_1,\ldots,v_r \in \sk$ we have
\[ \sum_{\substack{\K \in \widetilde{E}_G(k,d)\\
v_i \,\,\mbox{\tiny ramified in } \, \rho_i \circ \pp \,\, \forall \, i}}1 =H\prod_{i=1}^r \left(\kappa-\# \ker \rho_i\right)  H_{_{\sum_0}}'c_k^{\kappa-1} d^{\kappa-2}q^{d+\kappa-1} \left\{1+O_{G,k}\left(C^rd^{-1}\right)\right\}, \]
where the constants $H,c_k,C$ are the same as in Theorem \ref{3:rho} and
\[H_{\Sigma_0}'=\prod_{i=1}^r|v_i|^{-1}\left(1+(\kappa-1)|v|^{-1}\right)^{-1}.\]
\end{thm}

\noindent {\bf Proof.} Define
\[\Ba:=\sum_{\substack{\K \in \widetilde{E}_G(k,d)\\
v_i \,\,\mbox{\tiny ramified in } \, \rho_i \circ \pp \,\, \forall \, i}}1,\]
where $\vec{\rho}=(\rho_1,\ldots,\rho_r)$ and $\vec{v}=(v_1,\ldots,v_r)$. Since the idea of the proof is very similar to that of Theorem \ref{3:rho}, we only sketch the main steps.

First by class field theory we can rewrite $\Ba$ as
\[\Ba=\sum_{\substack{\chi:J/k^* \to G\\
v_i \,\,\mbox{\tiny ramified in } \, \rho_i \circ \pp \,\, \forall \, i\\
C(\chi)=q^d}} 1. \]
Let $\Sigma_0:=\{v_1,\ldots,v_r\}$ and let $\Sigma' \subset \sk$ be a finite subset order say $c$ such that $o_{\Sigma'}$ has class number 1. Then define
\[\Sigma=\Sigma_0 \bigcup \Sigma', \quad r'=\#\Sigma \le r+c. \]
The ring $o_{\Sigma}$ also has class number 1, so using the isomorphism $J_{\Sigma}/o_{\Sigma}^* \to J/k^*$ we have
\[\Ba=\sum_{\substack{\chi:J_{\Sigma}/o_{\Sigma}^* \to G\\
v_i \,\,\mbox{\tiny ramified in } \, \rho_i \circ \pp \,\, \forall \, i\\
C(\chi)=q^d}} 1. \]
Following the same argument as in the proof of Theorem \ref{3:rho}, we define the twists
\begin{eqnarray*} \label{4:twi} \Bcb=
\sum_{\substack{\chi:J_{\Sigma} \to G \\
v_i \, \mbox{\tiny ramified in } \rho_i \circ \chi \, \forall \, i\\
C(\chi)=q^d}} \dot{\chi}(\vec{\epsilon}), \end{eqnarray*}
where $\vec{\epsilon}=(\epsilon_1,\ldots,\epsilon_t) \in \A=\prod_{i=1}^t o_{\Sigma}^*/o_{\Sigma}^{*n_i}$, and given a $\chi: J_{\Sigma} \to G$ with projection $\chi_i: J_{\Sigma} \to \Z/n_i\Z$, we define $\dot{\chi}(\vec{\epsilon})=\prod_{i=1}^t
\zeta_{n_i}^{\chi_i(\epsilon_i)}$. We also have the relation \begin{eqnarray} \label{4:epf} \Ba=\frac{1}{\#\A}
\sum_{\vec{\epsilon} \in \A}\Bcb. \end{eqnarray}
The generating function of $\Bcb$ is
\[\Ccb=\sum_{\substack{\chi:J_{\Sigma} \to G \\
v_i \, \mbox{\tiny ramified in } \rho_i \circ \chi \, \forall \, i}} \frac{\dot{\chi}(\vec{\epsilon})}{C(\chi)^s}. \]
It has the Euler products
\begin{eqnarray*} F_G(\vec{\rho},\vec{v},\vec{\epsilon};s)&=&
\prod_{v \notin \Sigma }\left(\sum_{\chi_v:o_v^* \to G} \frac{\dot{\chi}_v(\vec{\epsilon})}{Nv^{c(\chi_v)s}}\right) \prod_{v \in \Sigma \setminus \Sigma_0}\left(\sum_{\chi_v:k_v^* \to G} \frac{\dot{\chi}_v(\vec{\epsilon})}{Nv^{c(\chi_v)s}}\right) \times \\
&& \prod_{i=1 }^r\left(\sum_{\substack{\chi_{v_i}:k_{v_i}^* \to G\\
v_i \, \mbox{\tiny ramified in } \rho_i \circ \chi}}  \frac{\dot{\chi}_{v_i}(\vec{\epsilon})}{Nv_i^{c(\chi_{v_i})s}}
\right).\end{eqnarray*}
Notice that $v$ is ramified in $\rho \circ \chi$ if and only if $\chi(y_{v}) \notin \ker \rho$ where $y_v$ is a generator of the tame inertia group of $k_v$. Similar to the arguments in Theorem \ref{3:rho}, we have:

\begin{lemma} \label{4:fg} If $\vec{\epsilon}=\vec{1}$, then
\[F_G(\vec{\rho},\vec{v},\vec{1};s)=\kappa^{r'} \prod_{i=1}^r \left(\kappa-\# \ker \rho_i\right) \prod_{v \notin \Sigma_0} \left(1+\frac{\kappa-1}{|v|^{s}}\right) \prod_{i=1}^r \frac{1}{|v_i|^{s}}. \]
If $\vec{\epsilon} \ne \vec{1}$, then
\[F_G(\vec{\rho},\vec{v},\vec{\lambda},\vec{\epsilon};s)= h(s) \prod_{v \notin \Sigma_0} \left(1-\frac{1}{|v|^{s}}\right) ,  \]
where
\[|h(s)| \le \kappa^{r'} \prod_{i=1}^r \left(1+\frac{\kappa}{|v_i|^{\Re s}}\right). \]
\end{lemma}
From Lemma \ref{4:fg}, similar to the arguments in the proofs of Lemmas \ref{3:ad1} and \ref{3:ad2}, we can obtain
\[\Bcb \ll_{G,k,\eta} C^{r} q^{(1+\eta)d/2}, \quad \forall \, \vec{\epsilon} \ne \vec{1}, \]
and
\[\Bca=\kappa^{r'} \prod_{i=1}^r \left(\kappa-\# \ker \rho_i\right) H \cdot H_{_{\sum_0}}'c_k^{\kappa-1}d^{\kappa-2} q^{d+\kappa-1} \left\{1+O_{G,k}\left(C^rd^{-1}\right)\right\}, \]
where the constants $H,c_k,C$ are the same as in Theorem \ref{3:rho}, and $H_{\Sigma_0}'$ is given by
\[H_{\Sigma_0}'=\prod_{i=1}^r |v_i|^{-1} \left(1+\frac{\kappa-1}{|v|}\right)^{-1}. \]
From these estimates and using the relation (\ref{4:epf}) we find that
\[\Ba=\prod_{i=1}^r \left(\kappa-\# \ker \rho_i\right) H \cdot \binom{d+\kappa-2}{\kappa-2}  H_{_{\sum_0}}'c_k^{\kappa-1}q^{d+\kappa-1} \left\{1+O_{G,k}\left(C^rd^{-1}\right)\right\}. \]
This completes the proof of Theorem \ref{4:rho}. \quad $\square$

From Theorems \ref{3:rho} and \ref{4:rho} we obtain immediately that
\begin{cor} \label{5:rho} For any non-principal characters $\rho_1,\ldots,\rho_r \in \widehat{G}$ of order say $\tau_1,\ldots,\tau_r$ respectively and any distinct places $v_1,\ldots,v_r \in \sk$.

\noindent {\bf (i).} Let  $\lambda_1,\ldots,\lambda_r$ be any integers such that $\tau_i \nmid \lambda_i$ for some $i$, then
\[ \frac{1}{\#\widetilde{E}_G(k,d) } \sum_{\K \in \widetilde{E}_G(k,d) } \prod_{i=1}^r \rho_i \circ \pp\left(v_i\right)^{\lambda_i} \ll_{G,k,\eta} C^r q^{(-1+\eta)d)/{2}} .\]

\noindent {\bf (ii).}
\[\frac{1}{\#\widetilde{E}_G(k,d)}\sum_{\substack{\K \in \widetilde{E}_G(k,d)\\
v_i \,\,\mbox{\tiny unramified in } \, \rho_i \circ \pp \,\, \forall \, i}}1= 1+O_{G,k}\left(\sum_{i=1}^r|v_i|^{-1}\right) . \]

\noindent {\bf (iii).}
\[\frac{1}{\#\widetilde{E}_G(k,d)} \sum_{\substack{\K \in \widetilde{E}_G(k,d)\\
v_i \,\,\mbox{\tiny ramified in } \, \rho_i \circ \pp \,\, \forall \, i}}1 \ll_{G,k} \prod_{i=1}^r |v_i|^{-1}. \]
\end{cor}

\noindent {\bf Proof.} (i) is immediate from (i) of Theorem \ref{3:rho} and (\ref{3:egkd}). From (ii) of Theorem \ref{3:rho} we obtain
\[\frac{1}{\#\widetilde{E}_G(k,d)}\sum_{\substack{\K \in \widetilde{E}_G(k,d)\\
v_i \,\,\mbox{\tiny unramified in } \, \rho_i \circ \pp \,\, \forall \, i}}1= \prod_{i=1}^r \frac{1+(\#\ker \rho_i-1)|v_i|^{-1}}{
1+(\kappa-1)|v_i|^{-1}} \left\{1+O_{G,k}\left(C^r d^{-1}\right)\right\}. \]
It is easy to see that
\[\prod_{i=1}^r \frac{1+(\#\ker \rho_i-1)|v_i|^{-1}}{
1+(\kappa-1)|v_i|^{-1}} =1+O_{G,k}\left(\sum_{i=1}^r|v_i|^{-1}\right).\]
(iii) is also obvious from Theorem \ref{4:rho}, as
\[\frac{1}{\#\widetilde{E}_G(k,d)} \sum_{\substack{\K \in \widetilde{E}_G(k,d)\\
v_i \,\,\mbox{\tiny ramified in } \, \rho_i \circ \pp \,\, \forall \, i}}1 \ll_{G,k}H_{\Sigma_0}' \ll_{G,k} \prod_{i=1}^r |v_i|^{-1}.\]
This completes the proof of Corollary \ref{5:rho}. \quad $\square$

\section{Distribution of zeta zeros for $\widetilde{E}_G(k,d)$}


In this section we will study the distribution of zeta zeros for function fields $\K$ running over the set $\widetilde{E}_G(k,d)$ as $d \to \infty$, where
\[\widetilde{E}_G(k,d):= \left\{\K: \mathrm{Gal}(\K/k) \overset{\pp}{\longrightarrow} G \, \mbox{ injective}, \deg_k \mathrm{Cond}(\K)=d\right\}. \]
More precisely, we will prove Theorem \ref{main} for $\widetilde{E}_G(k,d)$ instead of the set $E_G(k,d)$. This is most convenient, because certain estimates of sums over $\widetilde{E}_G(k,d)$ are provided by Corollary \ref{5:rho}. In the final Section distribution results on $E_G(k,d)$ will be derived from it.

The ideas of the proof are similar to those in \cite{fai,buc3, xio}. In particular readers may readily recognise that Corollary \ref{5:rho} obtained in the previous section looks similar to \cite[Theorem 2]{xio}, which was applied in an essential way to prove a general result there. Here Corollary \ref{5:rho} plays the same role and the arguments are similar. However, it seems not easy to simplify the proof. For the sake of completeness, we shall produce a proof with enough details. We follow closely the presentation of \cite[Section 3]{xio}.

\subsection{Preparation}

For a symmetric interval $\I=[-\beta/2,\beta/2]$, $0<\beta<1/2$, let $I^{\pm}_l(x)$ be the two Beurling-Selberg polynomials of degree $l$ defined in Section 2 such that $I^{-}_l \le \mathbf{1}_{\I} \le I^{+}_l$. Since $d |\I|=d \beta \to \infty$ as $d \to \infty$, we can choose integers $l=l(d)$ in such a way that
\begin{eqnarray*} \label{3:K} \frac{d}{l} \to \infty,\, l \beta \to \infty,\, \mbox{ and } \frac{d}{l } \ll \left(\log l\beta \right)^{1/4}  \, \mbox{ as } d \to \infty.\end{eqnarray*}
For any non-principal character $\rho:G \to \CC^*$ and any $\K \in \widetilde{E}_G(k,d)$, let
$\{\theta_{\rho,\K,i}: 1 \le i \le m_{\rho,\K}\} $ be the set of angles of zeta zeros for the L-function $L(k,\rho \circ \pp,s)$. From the monotonicity of $I^{\pm}_l$, we have
\begin{eqnarray*} \label{3:NK} N^{-}_{\rho,l}(\K) \le N_{\rho,\I} (\K) \le N^{+}_{\rho,l}(\K), \end{eqnarray*}
where
\[N^{\pm}_{\rho,l}(\K) = \sum_{i=1}^{m_{\rho,\K}} I^{\pm}_l\left(\theta_{\rho,\K,i}\right), \quad N_{\rho,\I}(\K)=\{i: \theta_{\rho,\K,i} \in \I\}\,. \]
Let $I_l(x)=\sum_{n}c(n)e(nx):=I^{\pm}_l(x)$. We collect several properties of the Fourier coefficients $c(n)$ from Section 2 as follows:
\begin{itemize}
\item[(i)] $c(n)=0$ if $|n|>l$.

\item[(ii)] $c(0)=\beta+O\left(l^{-1}\right)$.

\item[(iii)] $|c(n)n| \ll 1$ for any $n \in \Z$.

\item[(iv)] $\sum_{v \in \sk} c(\deg v)^2 (\deg v)^2 |v|^{-1} =\frac{1}{2 \pi^2} \log (l \beta) +O(1)$.
\item[(v)] $\sum_{n} c(2n) \ll 1$.
\end{itemize}
Define $N_{\rho,l}(\K):=N^{\pm}_{\rho,l}(\K)$. Then
\begin{eqnarray*} N_{\rho,l}(\K)&=&\sum_{i=1}^{m_{\rho,\K}}I_{l} \left(\theta_{\rho,\K,i}\right)=\sum_{n \in \Z} c(n) \sum_{i=1}^{m_{\rho,\K}} e\left(n \theta_{\rho,\K,i}\right)\,.\end{eqnarray*}
From the explicit formulas in (\ref{2:exp}) we obtain
\begin{eqnarray*} N_{\rho,l}(\K) &=& c(0) m_{\rho,\K} - \sum_{0 \ne n \in \Z}c(n) q^{-|n|/2} \sum_{\substack{v \\ \deg v|n}} (\deg v) \rho \circ \pp(v)^{n/\deg v} \,.\end{eqnarray*}
Since $|c(n)n| \ll 1$, we have
\[ \sum_{|n| \le l} |c(n)| q^{-|n|/2} \ll 1\,.\]
Using the fact that $\left|\rho \circ \pp(v)\right| \le 1$, and $c(n) =c(-n) \in \R$, we derive
\[N_{\rho,l}(\K)= \beta m_{\rho,\K}- \sum_{1 \le n \le l}c(n) q^{-n/2} \sum_{\substack{v \in \sk\\ \deg v|n}} (\deg v) \left\{ \rho \circ \pp(v)^{\frac{n}{\deg v}}+\rho \circ \pp(v)^{-\frac{n}{\deg v}}\right\} +O\left(\frac{m_{\rho,\K}}{l}\right)\,. \]
We may rewrite it as
\[N_{\rho,l}(\K)= \beta \MM +S_{\rho,l}(\K) +O\left(\frac{\MM}{l}\right)\,, \]
where
\[S_{\rho,l}(\K)=- \sum_{\substack{v \in \sk\\
r \ge 1}}c(r \deg v) |v|^{-r/2} (\deg v) \left\{ \rho \circ \pp(v)^{r}+\rho \circ \pp(v)^{-r}\right\}. \]
Now it is easy to see that
\[S_{\rho,l}(\K) \ll \sum_{1 \le n \le l} q^{-n/2} \sum_{\substack{v \in \sk\\
\deg v \mid n}} 1 \ll \sum_{1 \le n \le l} q^{n/2} \ll q^{l/2},\]
and hence
\[\left|N_{\rho,l}(\K)-\beta \GG  \right| \ll \frac{\GG}{l}+q^{l/2}\,.\]
Noting that $|\I|=\beta$, $\GG \asymp d$ and taking $l \asymp \log_q d-\log_q \log d$, we have deduced that the zeros are uniformly distributed:
\begin{proposition} As $d \to \infty$, for each $\K \in \widetilde{E}_G(k,d)$, every fixed symmetric interval $\I=[-\beta/2,\beta/2]$ contains asymptotically $\GG |\I|$ angles $\theta_{\rho,\K,i}$ for each $\rho$. In fact
\[N_{\rho,\I} (\K)=\GG |\I|+ O\left(\frac{\GG}{\log \GG}\right)\,. \]
\end{proposition}

Denote
\begin{eqnarray} \label{3:TK} T_{\rho,l}(\K)=- \sum_{v \in \sk} c(\deg v) \deg v |v|^{-1/2} \left\{ \rho \circ \pp(v)+\rho \circ \pp(v)^{-1}\right\}, \end{eqnarray}
and
\begin{eqnarray} \label{3:triangle} \triangle_{\rho,l}(\K)=- \sum_{v \in \sk } c(2\deg v) \deg v |v|^{-1} \left\{ \rho \circ \pp(v)^{2}+\rho \circ \pp(v)^{-2}\right\}. \end{eqnarray}
Then
\begin{eqnarray*} S_{\rho,l}(\K)-T_{\rho,l}(\K)-\triangle_{\rho,l}(\K)&=&-\sum_{\substack{v \in \sk\\
r \ge 3}} c(r \deg v) |v|^{-r/2} (\deg v) \left\{ \rho \circ \pp(v)^r+\rho \circ \pp(v)^{-r}\right\}. \end{eqnarray*}
It is easy to see that this is bounded by
\begin{eqnarray*}
\ll \sum_{\substack{v \in \sk\\
r \ge 3}} q^{-r \deg v/2} \le \sum_{r \ge 3} \sum_{n \le l/3} q^{-rn/2} q^n \ll 1.  \end{eqnarray*}
Therefore
\begin{eqnarray} \label{3:id2} N_{\rho,l}(\K)-\GG\beta=T_{\rho,l}(\K)+\triangle_{\rho,l}(\K)+O\left(\frac{\GG}{l}\right)\,.
\end{eqnarray}
We denote by $\langle \bullet \rangle$ the mean value of any quantity defined on $\FC$, that is, let $\chi: \FC \to \CC$ be a map, then
\[\langle \chi \rangle:=\frac{1}{\# \FC} \sum_{\K \in \FC} \chi(\K)\,. \]
The goal is to compute for any fixed nonnegative integers $r_1,r_2,\ldots,r_t$ the moment $\left\langle \prod_{j=1}^t\left(N_{\rho_j,l}(\bullet)-m_{\rho_j,\bullet} \beta\right)^{r_j}\right\rangle$. For this purpose we need to compute various moments for $\triangle_{\rho_j,l}(\bullet)$ and $T_{\rho_j,l}(\bullet)$ first. To simplify notation, in what follows the implied constants in ``$O$'' and ``$\ll $'' may depend on $G,k$ and the integer $r=\sum_j r_j$.

\subsection{Moments of $\triangle_{\rho,l}(\K)$}

For each $\rho$ of order $\tau \ge 2$, we may consider two cases.

\noindent {\bf Case (1). $\tau \ge 3$.} Then for each fixed positive integer $r$, we have from (\ref{3:triangle}) \[\triangle_{\rho,l}(\K)^{2r} =\sum_{v_1, \ldots, v_{2r} \in \sk} \prod_{i=1}^{2r} c(2 \deg v_i) |v_i|^{-1} \deg v_i \sum_{\lambda_1, \ldots, \lambda_{2r} \in \left\{1,-1\right\}} \prod_{i=1}^{2r} \rho \circ \pp(v_i)^{2\lambda_i}\,. \]
Hence
\[\langle \left(\triangle_{\rho,l}\right)^{2r} \rangle = \sum_{v_1, \ldots, v_{2r} \in \sk} \prod_{i=1}^{2r} c(2 \deg v_i) |v_i|^{-1} \deg v_i \sum_{\lambda_1, \ldots, \lambda_{2r} \in \left\{1,-1\right\}}  \left\langle \prod_{i=1}^{2r} \rho \circ \pp(v_i)^{2 \lambda_i}\right\rangle \,. \]
Consider $v_1^{2 \lambda_1} \cdots v_{2r}^{2 \lambda_{2r}}$ as an element of $\DD$, the free abelian multiplicative group generated by all the places of $k$. If $v_1^{2\lambda_1} \cdots v_{2r}^{2\lambda_{2r}}$ is not a $\tau$-th power in $\DD$, then by (i) of Corollary \ref{5:rho}, we have
\[\left\langle \prod_{i=1}^{2r} \rho \circ \pp(v_i)^{2 \lambda_i}\right\rangle \ll_{\eta} C^r q^{(-1+\eta)d/2}. \]
Since $d/l \to \infty$, the total contribution in this case is
\[\langle (\triangle_{\rho,l})^{2r} \rangle_1 \ll_{\eta} \sum_{\substack{v_1, \ldots, v_{2r}\\
\deg v_i \le l}} \prod_{i=1}^{2r} |v_i|^{-1} \sum_{\lambda_1,\ldots,\lambda_{2r} \in \{1,-1\}} C^{2r} q^{(-1+\eta)d/2} \ll 1. \]
If $v_1^{2\lambda_1} \cdots v_{2r}^{2\lambda_{2r}}$ is a $\tau$-th power in $\DD$, since $\tau \ge 3$, for this to happen, each $v_i$ must be paired off with at least one $v_j, i \ne j$, $v_i=v_j$. Hence the total contribution in this case is at most
\[\langle \left(\triangle_{\rho,l}\right)^{2r} \rangle_2 \ll \left(\sum_{v} |v|^{-2}\right)^{r} \ll 1\,. \]

\noindent {\bf Case (2) $\tau=2$}. Then
\[\triangle_{\rho,l}(\K) =-2 \sum_{\substack{v \in \sk\\
v \, \mbox{\tiny unramified in } \rho \circ \pp}} c(2 \deg v) |v|^{-1} (\deg v) =\triangle_{1}(\K)+\triangle_{2}(\K), \]
where
\begin{eqnarray*} \triangle_{1}(\K)&=&-2 \sum_{\substack{v \in \sk}} c(2 \deg v) |v|^{-1} (\deg v)
=-2 \sum_{1 \le n \le l/2}c(2n) n q^{-n} \sum_{\substack{v \in \sk\\
\deg v =n}}1 \\
&=& -2 \sum_{1 \le n \le l/2}c(2n) n q^{-n} \left(q^n/n+O(q^{n/2})\right) \\
&=& -2 \sum_{1 \le n \le l/2} c(2n)+O(1) \ll 1, \mbox{ from (v)}, \end{eqnarray*}
and
\begin{eqnarray*} \triangle_{2}(\K)&=&2 \sum_{\substack{v \in \sk\\
v \, \mbox{\tiny ramified in } \rho \circ \pp}} c(2 \deg v) |v|^{-1} (\deg v).
\end{eqnarray*}
Since $|c(n)n| \ll 1$ we find
\begin{eqnarray*} \left \langle \left|(\triangle_{2})^r\right| \right \rangle & \ll & \left \langle \sum_{\substack{v_1,\ldots,v_r \in \sk\\
v_i \, \mbox{\tiny ramified in } \rho \circ \pp \, \forall i}} |v_1|^{-1} \cdots |v_r|^{-1} \right \rangle \\
&=& \sum_{\substack{v_1,\ldots,v_r \in \sk
}} |v_1|^{-1} \cdots |v_r|^{-1} \frac{1}{\# \widetilde{E}_G(k,d)} \sum_{\substack{\K \in \widetilde{E}_G(k,d)\\
v_i \, \mbox{\tiny ramified in } \rho \circ \pp \, \forall i}} 1.
\end{eqnarray*}
Now by (iii) of Corollary \ref{5:rho} we obtain
\begin{eqnarray*} \left \langle \left|(\triangle_{2})^r\right| \right \rangle & \ll & \sum_{\substack{v_1,\ldots,v_r \in \sk
}} |v_1|^{-1} \cdots |v_r|^{-1} \prod_{v \in \{v_1,\ldots,v_r\}} |v|^{-1} \\
& \ll & \left(\sum_{v \in \sk} |v|^{-2}\right)^r \ll 1.
\end{eqnarray*}
Combining the two cases above we conclude that
\begin{eqnarray} \label{3:tri} \langle (\triangle_{\rho,l})^{2r} \rangle \ll 1\,. \end{eqnarray}

\subsection{Moments of $T_{\rho,l}(\K)$}

For each fixed positive integer $r$, from (\ref{3:TK}) we have
\[ T_{\rho,l}(\K)^r =(-1)^r \sum_{v_1, \ldots, v_r \in \sk} \prod_{i=1}^r c(\deg v_i) |v_i|^{-1/2} (\deg v_i) \sum_{\lambda_1, \ldots, \lambda_r \in \{1,-1\} } \rho \circ \pp \left(v_1^{\lambda_1} \cdots v_{r}^{\lambda_{r}}\right)\,. \]
Hence
\[\langle \left(T_{\rho,l}\right)^{r} \rangle = (-1)^r \sum_{v_1, \ldots, v_r \in \sk} \prod_{i=1}^r c(\deg v_i) |v_i|^{-1/2} (\deg v_i) \sum_{\lambda_1, \ldots, \lambda_r \in \{1,-1\} } \left\langle \rho \circ \pp \left(v_1^{\lambda_1} \cdots v_{r}^{\lambda_{r}}\right) \right \rangle \,. \]
If $v_1^{\lambda_1} \cdots v_{r}^{\lambda_{r}}$ is not a $\tau$-th power in $\DD$, the total contribution in this case, from (i) of Corollary \ref{5:rho}, is
\[\langle \left(T_{\rho,l} \right)^{r} \rangle_1 \ll \sum_{\substack{v_1, \ldots, v_{r} \in \sk\\
\deg v_i \le l}} \prod_{i=1}^{r} |v_i|^{-1/2} \sum_{\lambda_1, \ldots, \lambda_{r} \in \{1,-1\}} C^r q^{(-1+\eta)d/2} \ll 1. \]
If $v_1^{\lambda_1} \cdots v_{r}^{\lambda_{r}}$ is a $\tau$-th power in $\DD$, for this to happen, each $v_i$ must be paired off with another $v_j, i \ne j$, $v_i = v_j$. The most ``economic'' way of doing that is that each $v_i$ is paired off with exactly one $v_j, i \ne j, v_i=v_j$, so that $r=2s$ must be even. The number of choices for such arrangement is $(2s)!/s! 2^s$. Moreover, the exponents must satisfy the condition $\lambda_i+\lambda_j \equiv 0 \pmod{\tau}$, and there are exactly $2 r_{\rho}$ choices for $\lambda_i, \lambda_j \in \{1,-1\}$, where $r_{\rho}=1$ if $\tau \ge 3$, and $r_{\rho}=2$ if $\tau=2$, where $\tau$ is the order of $\rho$. Hence the contribution in this case is
\[\langle \left(T_{\rho,l}\right)^{r} \rangle_0 = (2 r_{\rho})^s\frac{(2s)!}{s! 2^s}\sum_{\substack{v_1, \ldots, v_{s} \in \sk \\
\mbox{\tiny distinct } }} \prod_{i=1}^{s} c(\deg v_i)^2 |v_i|^{-1} (\deg v_i)^2 \left\langle \rho^{\tau} \circ \pp (v_1 \cdots v_s) \right \rangle \,. \]
Since $\rho$ has order $\tau$, from (ii) of Corollary \ref{5:rho} we find that
\[\left\langle\rho^{\tau} \circ \pp (v_1 \cdots v_s) \right\rangle=\frac{1}{\#\widetilde{E}_G(k,d)} \sum_{\substack{\K \in \widetilde{E}_G(k,d)\\
v_i \,\,\mbox{\tiny unramified in } \, \rho_i \circ \pp \,\, \forall \, i}}1= 1+O\left(\sum_{i=1}^s |v_i|^{-1}\right).\]
The contribution from the error term $O\left(\sum_{i=1}^s |v_i|^{-1}\right)$ is bounded by
\begin{eqnarray*} E \ll \left(\sum_{v \in \sk} c(\deg v)^2 (\deg v)^2 |v|^{-1}\right)^{s-1} \left(\sum_{v \in \sk} |v|^{-2} \right). \end{eqnarray*}
Using property (iv) of $c(n)$ and the estimate $\sum_{v \in \sk} |v|^{-2} \ll 1$, we find
\[E \ll \left(\log l \beta\right)^{s-1}\,.\]
The main term is
\[\frac{(r_{\rho})^s(2s)!}{s! }\sum_{\substack{v_1, \ldots, v_{s} \in \sk
\\ \mbox{\tiny distinct } }} \prod_{i=1}^{s} c(\deg v_i)^2 |v_i|^{-1} (\deg v_i)^2  \,.\]
Now we remove the restriction that $v_1, \ldots, v_s$ are distinct, introducing again an error of $O\left(\left(\log l \beta\right)^{s-2}\right)$. This gives us
\[\langle \left(T_{\rho,l}\right)^{r} \rangle_0= \frac{(r_{\rho})^s(2s)!}{s! } \left(\sum_{v \in \sk} c(\deg v)^2 (\deg v)^2 |v|^{-1}\right)^{s} +O\left(\left(\log l \beta\right)^{s-1}\right).\]
Using property (iv) of $c(n)$ again we derive that
\begin{eqnarray*} \langle \left(T_{\rho,l}\right)^{r} \rangle_0&=& \frac{(r_{\rho})^s(2s)!}{s!} \left(\frac{1}{2 \pi^2} \log (l \beta)+O(1)\right)^{s} +O\left(\left(\log l \beta\right)^{s-1}\right), \end{eqnarray*}
and from it we obtain
\begin{eqnarray*} \label{3:main}
\langle \left(T_{\rho,l}\right)^{r} \rangle_0 & = &  \frac{ (r_{\rho})^s(2s)!}{2^s \pi^{2s} s! } \left(\log l\beta\right)^s + O\left(\left(\log l \beta\right)^{s-1}\right),\end{eqnarray*}
where $r=2s$ is an even number.

On the other hand, if each $v_i$ is paired off with another $v_j, i \ne j$, $v_i=v_j$, but it is not the most ``economic'', that is, there is one $v_i$ which is paired off with at least two others $v_{j_1}, v_{j_2}, i \ne j_1 \ne j_2$ with $v_i=v_{j_1}=v_{j_2}$. Similar to \cite[4.2.2 Case two]{xio}, we find easily that the total contribution from this case is bounded by
\[\ll (\log l \beta)^{(r-3)/2}. \]
We conclude that
\[\left\langle (T_{\rho,l})^r \right\rangle =\frac{(r_{\rho})^{r/2}\delta(r) r!}{2^{r/2}\pi^r (r/2)!}(\log l \beta)^{r/2} +O\left((\log l \beta)^{-1+r/2}\right),\]
where $\delta(r)=1$ if $r$ is even and $\delta(r)=0$ if $r$ is odd.

\subsection{General moments of $T_{\rho,l}$}
Let $\rho_1,\ldots,\rho_t \in \widehat{G}$ be non-principal characters of order $\tau_1,\ldots,\tau_t$ respectively. For any nonnegative integers $r_1, \ldots, r_t$, let $r=\sum_{j=1}^t r_j$. We have
\[\prod_{j=1}^t T_{\rho_j,l}(\K)^{r_j}=(-1)^{r}\sum_{v_{j,i}} \prod_{j,i}c\left(\deg v_{j,i}\right)|v_{j,i}|^{-1/2}\deg v_{j,i} \sum_{\lambda_{j,i} \in \{1,-1\}}\prod_{j,i} \rho_j \circ \pp\left(v_{j,i}^{\lambda_{j,i}} \right)\,,\]
where the index $j,i$ run over the range $1 \le j \le t$ and $1 \le i \le r_j$. Hence
\[\left\langle \prod_{j=1}^t \left(T_{\rho_j,l}\right)^{r_j} \right\rangle = (-1)^r \sum_{v_{j,i}} \prod_{j,i}c\left(\deg v_{j,i}\right)|v_{j,i}|^{-1/2}(\deg v_{j,i}) \sum_{\lambda_{j,i} \in \{1,-1\}}\left\langle \prod_{j,i} \rho_j \circ \pp\left(v_{j,i}^{\lambda_{j,i}} \right) \right \rangle \,. \]

Similarly, if some $v_{j,i}$ is not paired off with another $v_{j',i'}$, then by (i) of Corollary \ref{5:rho}, the total contribution in this case is bounded by $O(1)$; if each $v_{j,i}$ is paired off with another $v_{j',i'}$ and there is one element that is paired off with at least two other elements, using similar argument the total contribution in this case is bounded by $O\left(\left(\log K\beta\right)^{(r-3)/2}\right)$. The remaining cases are that each $v_{j,i}$ is paired off with exactly one $v_{j',i'}$ where $(j,i) \ne (j',i')$. Since $\rho_i \rho_j \ne 1$ for any $i \ne j$, if there are $(j,i),(j',i')$ with $j \ne j'$ such that $v_{j,i}=v_{j',i'}$, then by (i) of Corollary \ref{5:rho}, the total contribution again is bounded by $O(1)$. The main contribution comes from the case that for each $j$, $v_{j,i}$ is paired off with exactly one $v_{j,i'}$ with $i \ne i'$. For this to happen, then for each $j$, $r_j=2s_j$ must be even, and the number of choices of $v_{j,i}$ and $\lambda_{j,i}$ is $(2 r_{\rho_j})^{s_j}\frac{(2s_j)!}{2^{s_j} (s_j)!}$. Hence the total contribution in this case is
\[\left\langle \prod_{j=1}^{t}\left(T_{\rho,l}\right)^{r_j} \right\rangle_0 = \prod_{j=1}^t (2 r_{\rho_j})^{s_j}\frac{(2s_j)!}{(s_j)! 2^{s_j}}\sum_{\substack{v_{j,i} \in \sk
\\ \mbox{\tiny all distinct } }} \prod_{j,i} c(\deg v_{j,i})^2 |v_{j,i}|^{-1} (\deg v_{j,i})^2 \left\langle \prod_{j,i} \rho_j^{\tau_j} \circ \pp \left(v_{j,i} \right) \right \rangle  \,, \]
where the index $j,i$ runs over the range $1 \le j \le t$ and $1 \le i \le s_j$. Since
\[\left\langle \prod_{j,i} \rho_j^{\tau_j} \circ \pp \left(v_{j,i} \right) \right \rangle =1+ O\left(\sum_{j,i} |v_{j,i}|^{-1}\right)\,,  \]
the error term arising from $O\left(\sum_{j,i}|v_{j,i}|^{-1}\right)$ is $E \ll \left(\log K \beta\right)^{-1+\sum_{j}{s_j}}$, and the main term is
\[\prod_{j=1}^t (2r_{\rho_j})^{s_j}\frac{(2s_j)!}{(s_j)! 2^{s_j}}\sum_{\substack{v_{j,i} \in \sk \\
\mbox{\tiny all distinct } }} \prod_{j,i} c(\deg v_{j,i})^2 |v_{j,i}|^{-1} (\deg v_{j,i})^2 \,. \]
We may remove the restriction that all $v_{j,i}$'s are distinct, introducing again an error of $O\left(\left(\log K \beta\right)^{-2+\sum_j s_j}\right)$. This gives us
\[\left\langle \prod_{j=1}^t\left(T_{\rho_j,l}\right)^{r_j} \right\rangle_0= \prod_{j=1}^t\frac{(r_{\rho_j})^{s_j}(2s_j)!}{(s_j)! } \left(\sum_{v \in \sk} c(\deg v)^2 (\deg v)^2 |v|^{-1}\right)^{s_j} +O\left(\left(\log l \beta\right)^{-1+\sum_j s_j}\right).\]
From it we obtain
\begin{eqnarray*}
\left\langle \prod_{j=1}^t \left(T_{\rho_j,l}\right)^{r_j} \right\rangle_0 & = &  \prod_{j=1}^t\frac{ (r_{\rho_j})^{s_j}(2s_j)!}{2^{s_j} \pi^{2s_j} s_j! } \left(\log l \beta\right)^{s_j} + O\left(\left(\log l \beta\right)^{-1+\sum_j s_j}\right).\end{eqnarray*}
Combining the above estimates together we conclude that
\begin{eqnarray} \label{3:main2}
\left\langle \prod_{j=1}^{t} \left(T_{\rho_j,l}\right)^{r_j} \right\rangle & = &  \prod_{j=1}^{t}\frac{ (r_{\rho_j})^{r_j/2}\delta(r_j) (r_j)!}{2^{r_j/2} \pi^{r_j} (r_j/2)! } \left(\log l \beta\right)^{r_j/2} + O\left(\left(\log l \beta\right)^{-1+r/2}\right),\end{eqnarray}
where $\delta(s)=1$ if $s$ is even and $\delta(s)=0$ if $s$ is odd, and $r=\sum_{j=1}^{t} r_j$, and the implied constant in ``$O$'' may depend on $G,k$ and $r$.

\subsection{Proof of Theorem \ref{main} for $\widetilde{E}_G(k,d)$}
Very similar to \cite[Section 5]{xio}, combining the results obtained in (\ref{3:main2}), (\ref{3:tri}) and using (\ref{3:id2}) we can obtain
\begin{eqnarray} \label{3:fin} \left\langle \prod_{j=1}^t \left(\frac{N_{\rho_j,l}(\bullet)-\beta m_{\rho,\bullet}}{\sqrt{ \frac{r_{\rho_j}}{\pi^2}\log (l\beta)}}\right)^{r_j} \right\rangle= \prod_{j=1}^t \frac{\delta(r_j) r_j! }{2^{r_j/2} \left(r_j/2\right)! } +O\left(\left(\log l\beta\right)^{-1/4} \right)\,. \end{eqnarray}
Since $N_{\rho_j,l}^{-}(\K) \le N_{\rho_j,\I}(\K) \le N_{\rho_j,l}^{+}(\K)$ for each $\rho_j$, we have taken for granted in \cite[Section 5]{xio} that in Equation (\ref{3:fin}) the terms $N_{\rho_j,l}(\bullet)=N_{\rho_j,l}^{\pm}(\bullet)$ could be replaced by $N_{\rho_j,\I}(\bullet)$ and Equation (\ref{3:fin}) still holds true. Since for a standard Gaussian distribution, the odd moments vanish and the even moments are
\[\frac{1}{\sqrt{2 \pi}} \int_{\infty}^{\infty} x^{2r}e^{-x^2/2} \ud \, x = \frac{(2r)!}{2^r r!}\,,\]
we shall conclude directly that as $\K$ runs through $\widetilde{E}_G(k,d)$ with $d \to \infty$,
\[ \left(\frac{N_{\rho_1,\I}(\K)-m_{\rho_1,\K} \beta}{\sqrt{\frac{r_{\rho_1}}{\pi^2} \log m_{\rho_1,\K} \beta}}, \ldots, \frac{N_{\rho_t,\I}(\K)-m_{\rho_t,\K} \beta}{\sqrt{\frac{r_{\rho_t}}{\pi^2} \log m_{\rho_t,\K} \beta}} \right)\]
converges weakly to $t$ identical independent standard Gaussian distributions.

That in Equation (\ref{3:fin}) the terms $N_{\rho_j,l}(\bullet)=N_{\rho_j,l}^{\pm}(\bullet)$ could be replaced by $N_{\rho_j,\I}(\bullet)$ and Equation (\ref{3:fin}) still holds true -- rigorously speaking, that is correct but not quite obvious and needs to be proved. We take this opportunity to provide an argument here.

Define
\begin{eqnarray*} W_{\rho,l}(\K) &:=&N_{\rho,l}^{+}(\K)-N_{\rho,l}^{-}(\K).\end{eqnarray*}
We have from (\ref{3:id2})
\begin{eqnarray*}
W_{\rho,l}(\K)&=& W'_{\rho,l}(\K)+W''_{\rho,l}(\K) + O\left(\frac{d}{l}\right),
\end{eqnarray*}
where
\[W'_{\rho,l}(\K)=\sum_{v \in \sk} \frac{(\deg v) \left(c^-(\deg v)-c^{+}(\deg v)\right)}{|v|^{1/2}}\left\{\rho \circ \pp(v)+\rho \circ \pp(v)^{-1}\right\},\]
and
\[W''_{\rho,l}(\K)=\sum_{v \in \sk} \frac{(\deg v) \left(c^-(2\deg v)-c^{+}(2\deg v)\right)}{|v|}\left\{\rho \circ \pp(v)^2+\rho \circ \pp(v)^{-2}\right\}.\]
Since
\[\left|c^+(n)-c^-(n)\right| \le \frac{1}{2(l+1)}, \]
we find that
\[W''_{\rho,l}(\K) \ll \frac{1}{l} \sum_{\substack{v \in \sk\\
\deg v \le l}}(\deg v)|v|^{-1}\ll \frac{1}{l} \sum_{n \le l}n q^{-n} \frac{q^n}{n} \ll 1. \]
Hence
\begin{eqnarray} \label{6:wrho} W_{\rho,l}(\K)=W'_{\rho,l}(\K)+O\left(\frac{d}{l}\right). \end{eqnarray}
Next for any non-negative integers $r_1,\ldots,r_t$, let $r=\sum_{j} r_j$, we obtain
\[\left\langle \prod_{j=1}^t(W'_{\rho_j,l})^{2r_j} \right\rangle  \ll l^{-2r} \sum_{\substack{v_{j,i} \in \sk\\
\deg v_{j,i} \le l}} \prod_{j,i} \frac{\deg v_{j,i}}{|v_{j,i}|^{1/2}} \sum_{\substack{\lambda_{j,i} \in \{1,-1\}}} \left\langle \prod_{j,i} \rho_j \circ \pp(v_{j,i})^{\lambda_{j,i}} \right\rangle,\]
where the index $j,i$ runs over the range $1 \le j \le t$ and $1 \le i \le 2r_j$.

Again, if one $v_{j,i}$ is not paired off with any other $v_{j',i'}$, then by (i) of Corollary \ref{5:rho}, the contribution for this case is bounded by $O(1)$. If for any $j,i$, the place $v_{j,i}$ is paired off with another $v_{j',i'}$, then the total contribution for this case is bounded by
\begin{eqnarray*} & \ll & l^{-2r} \left(\sum_{\substack{v \in \sk\\
\deg v \le l}}(\deg v)^2 |v|^{-1}\right)^r = l^{-2r} \left(\sum_{1 \le n \le l}n^2 q^{-n}\sum_{\substack{v \in \sk\\
\deg v=n}}1 \right)^r \\
& \ll & l^{-2r} \left(\sum_{1 \le n \le l}n^2 q^{-n} \frac{q^n}{n} \right)^r \ll 1. \end{eqnarray*}
So we have
\[\left\langle \prod_{j=1}^t(W'_{\rho_j,l})^{2r_j} \right\rangle \ll 1.\]
Using (\ref{6:wrho}) we obtain
\begin{eqnarray*} \left|\left\langle \prod_{j=1}^t(W_{\rho_j,l})^{r_j} \right\rangle \right|&\ll & \sum_{\substack{u_j+v_j=r_j , \, \forall j}} \left\langle \left|\prod_{j=1}^t(W'_{\rho_j,l})^{u_j} \left(\frac{d}{l}\right)^{v_j}\right| \right\rangle \\
& \ll & \sum_{\substack{u_j+v_j=r_j , \, \forall j}} \left(\frac{d}{l}\right)^{\sum_j v_j} \left\langle \prod_{j=1}^t(W'_{\rho_j,l})^{2u_j} \right\rangle^{1/2}.
\end{eqnarray*}
Since $d/l \ll (\log l \beta)^{1/4}$, this gives
\begin{eqnarray} \label{6:wrho2} \left|\left\langle \prod_{j=1}^t(W_{\rho_j,l})^{r_j} \right\rangle \right|
& \ll & \sum_{\substack{u_j+v_j=r_j , \, \forall j}} \left(\frac{d}{l}\right)^{\sum_j v_j} \ll (\log l \beta)^{r/4}.
\end{eqnarray}
With this estimate, and noting that
\[N_{\rho_j,\I}(\K)- \beta m_{\rho_j,\K} =N_{\rho_j,l}^{-}(\K)-\beta m_{\rho_j,\K}+V_{\rho_j,l}(\K),\]
where
\[0 \le V_{\rho_j,l}(\K) = N_{\rho_j,\I}(\K)-N_{\rho_j,l}^{-}(\K) \le W_{\rho_j,l}(\K), \]
we find that
\begin{eqnarray*} \label{6:pn} \prod_{j=1}^t\left(N_{\rho_j,\I}(\K)- \beta m_{\rho_j,\K}\right)^{r_j} =
\prod_{j=1}^t \left(N_{\rho_j,l}^-(\K)- \beta m_{\rho_j,\K}\right)^{r_j} +E(\K),\end{eqnarray*}
where
\begin{eqnarray*} |E(\K)| &\ll & \sum_{\substack{u_j+v_j=r_j, \, \forall j \\
\sum_j v_j \ge 1}} \left|\prod_{j=1}^t \left(N_{\rho,l}^-(\K)- \beta m_{\rho,\K}\right)^{u_j}  W_{\rho_j,l}(\K)^{v_j}\right|. \end{eqnarray*}
Using the Cauchy-Schwartz inequality, we obtain
\begin{eqnarray*} \left\langle |E(\bullet)|\right\rangle
&\ll & \sum_{\substack{u_j+v_j=r_j, \, \forall j \\
\sum_j v_j \ge 1}} \left\langle \prod_{j=1}^t \left(N_{\rho,l}^-(\bullet)- \beta m_{\rho,\bullet}\right)^{2u_j} \right\rangle^{1/2} \left\langle \prod_{j=1}^t (W_{\rho_j,l})^{2v_j} \right\rangle^{1/2}.\end{eqnarray*}
From (\ref{3:fin}) and (\ref{6:wrho2}) we have
\[\left\langle |E(\bullet)|\right\rangle \ll \sum_{\substack{u_j+v_j=r_j, \, \forall j \\
\sum_j v_j \ge 1}} (\log l \beta)^{\sum_j u_j/2} (\log l \beta)^{\sum_{j}v_j/4} \ll (\log l \beta)^{\frac{r}{2}-\frac{1}{4}}. \]
From the above and (\ref{3:fin}) we conclude that
\begin{eqnarray} \label{4:fin} \left\langle \prod_{j=1}^t \left(\frac{N_{\rho_j,\I}(\bullet)-\beta m_{\rho_j,\bullet}}{\sqrt{ \frac{r_{\rho_j}}{\pi^2}\log (l\beta)}}\right)^{r_j} \right\rangle= \prod_{j=1}^t\frac{\delta(r_j) (r_j)! }{2^{r_j/2} \left(r_j/2\right)! } +O\left(\left(\log l\beta\right)^{-1/4} \right)\,, \end{eqnarray}
where the implied constant in ``$O$'' may depend on $G,k$ and $r=\sum_{j}r_j$. In the denominators we can also replace $l$ by $m_{\rho_j,\K}$ because
\[\log (l \beta) \sim \log (m_{\rho_j,\K} \beta). \]
This concludes the proof of Theorem \ref{main} for $\widetilde{E}_G(k,d)$ as $d \to \infty$. \quad $\square$

\section{Proof of Theorems \ref{main2} and \ref{main}}
Theorem \ref{main} on $E_G(k,d)$ can be derived easily from the result on $\widetilde{E}_G(k,d)$. For simplicity we define
\[a_{\rho}(\K):=\frac{N_{\rho,\I}(\K)-\beta m_{\rho,\K}}{\sqrt{\frac{r_{\rho}}{\pi^2}\log (l \beta)}}. \]
Equation (\ref{4:fin}) shows that
\[\frac{1}{\#\widetilde{E}_G(k,d)}\sum_{\substack{\K \in \widetilde{E}_G(k,d)}}\prod_{j=1}^t a_{\rho_j}(\K)^{r_j}=\prod_{j=1}^t\frac{\delta(r_j) (r_j)! }{2^{r_j/2} \left(r_j/2\right)! }+O\left((\log l \beta)^{-1/4}\right). \]
We need to prove that the above asymptotic formula holds true if we replace $\widetilde{E}_G(k,d)$ by $E_G(k,d)$. First, it is clear by definition that
\[E'_G(k,d):=\widetilde{E}_G(k,d) \setminus E_G(k,d)=\bigcup_{H \, \mbox{\tiny proper subgroup of} \, G} \widetilde{E}_H(k,d). \]
From the asymptotic formula (\ref{3:egkd})
\[\#E'_G(k,d) \le \sum_{H \, \mbox{\tiny proper subgroup of} \, G} \# \widetilde{E}_H(k,d) \ll \frac{1}{d} \cdot \#\widetilde{E}_G(k,d),\]
hence we have
\[\#E_{G}(k,d)=\widetilde{E}_G(k,d) \left(1+O(d^{-1})\right). \]
Therefore
\[\frac{1}{\#{E}_G(k,d)}\sum_{\substack{\K \in \widetilde{E}_G(k,d)}}\prod_{j=1}^t a_{\rho_j}(\K)^{r_j}=\prod_{j=1}^t\frac{\delta(r_j) (r_j)! }{2^{r_j/2} \left(r_j/2\right)! }+O\left((\log l \beta)^{-1/4}\right). \]
Now we estimate
\[B=\frac{1}{\#{E}_G(k,d)}\sum_{\substack{\K \in E'_G(k,d)}}\prod_{j=1}^t a_{\rho_j}(\K)^{r_j}.\]
We find from (\ref{4:fin}) that
\begin{eqnarray*} |B| &\le & \sum_{\substack{H \, \mbox{\tiny proper subgroup of }\, G}} \left\{\frac{\# \widetilde{E}_H(k,d)}{\#E_G(k,d)}\right\} \frac{1}{ \# \widetilde{E}_H(k,d)} \sum_{\K \in \widetilde{E}_H(k,d)} \prod_{j=1}^t a_{\rho_j}(\K)^{r_j}\\
& \ll & d^{-1}. \end{eqnarray*}
Thus we obtain
\[\frac{1}{\#E_G(k,d)}\sum_{\substack{\K \in E_G(k,d)}} \prod_{j=1}^t a_{\rho_j}(\K)^{r_j}=\prod_{j=1}^t\frac{\delta(r_j) (r_j)! }{2^{r_j/2} \left(r_j/2\right)! }+O\left((\log l \beta)^{-1/4}\right). \]
This shows that the value $\left(a_{\rho_j}(\K)\right)_{j=1}^t$ as $\K$ runs through $E_G(k,d)$ with $d \to \infty$ converges weakly to $t$ independent standard Gaussian distributions. The completes the proof of Theorem \ref{main}. \quad $\square$

Theorem \ref{main2} can be derived from Theorem \ref{main} as follows: first notice that
\begin{eqnarray} \label{7:ni} N_{\I}(\K)=\sum_{1 \ne \rho \in \widehat{G}} N_{\rho,\I}(\K)+O(1),\end{eqnarray}
where the term $O(1)$ comes from zeros of $\zeta_k(s)$. We choose a subset $A \subset \widehat{G} \setminus \{1\}$ with the property that for any $1 \ne \rho \in \widehat{G}$, then
\[\#\left(\left\{\rho,\rho^{-1}\right\} \bigcap A \right)=1. \]
Obviously $A$ exists and there are many choices for $A$. Then (\ref{7:ni}) can be written as
\[N_{\I}(\K)=\sum_{\rho \in A} \epsilon_{\rho} N_{\rho,\I}(\K)+O(1),\]
where
\[ \epsilon_{\rho}=\left\{\begin{array}{lll}1&:& \mbox{ if } \rho^2=1;\\
2&:& \mbox{ if } \rho^2 \ne 1. \end{array}\right.\]
This is because if $\rho^2 \ne 1$, then $\rho \ne \rho^{-1}$ and $N_{\rho,\I}(\K)=N_{\rho^{-1},\I}(\K)$ for symmetry. We also have
\[N_{\I}(\K)-g_{\K} \beta=\sum_{\rho \in A} \epsilon_{\rho} \left(N_{\rho,\I}(\K)-m_{\rho,\K} \beta\right)+O(1). \]
From Theorem \ref{main}, as $\K$ runs through $E_G(k,d)$ with $d \to \infty$, the value $N_{\rho,\I}(\K)-m_{\rho,\K} \beta$ converges weakly to independent Gaussian distributions with mean zero and variance $\frac{r_{\rho}}{\pi^2} \log (l \beta)$ for all $\rho \in A$. So $N_{\I}(\K)-g_{\K} \beta$ also converges to a Gaussian distribution with mean zero. The variance is
\begin{eqnarray*} \sigma^2 &\sim & \sum_{\rho \in A} \epsilon_{\rho}^2 \cdot \frac{r_{\rho}}{\pi^2} \log (l \beta) =\sum_{\rho \in A, \rho^2=1} \frac{2}{\pi^2} \log (l \beta)+\sum_{\rho \in A, \rho^2 \ne 1} \frac{4}{\pi^2} \log (l \beta). \end{eqnarray*}
It is easy to see that the right hand side is
\begin{eqnarray*}
\sum_{\rho \in \widehat{G}\setminus \{1\}, \rho^2 \ne 1} \frac{2}{\pi^2} \log (l \beta) +\sum_{\rho \in \widehat{G}, \rho^2 \ne 1} \frac{2}{\pi^2} \log (l \beta)
&=& \frac{2(\kappa-1)}{\pi^2} \log (l \beta),
\end{eqnarray*}
where $\kappa=\#G$. We can also replace $l$ by $g_{\K}$ since
\[\log (l \beta) \sim \log (g_{\K} \beta). \]
This completes the proof of Theorem \ref{main2}. \quad $\square$



\begin{thebibliography}{99}



\bibitem{buc1} A. Bucur, C. David, B. Feigon, M. Lal\'in, \emph{Statistics for traces of cyclic trigonal curves over finite fields}, Int. Math. Res. Not. IMRN (2010), No. 5, 932--967.

\bibitem{buc11} A. Bucur, C. David, B. Feigon, M. Lal\'in, \emph{Fluctuations in the number of points on smooth plane curves over finite fields}, J. Number Theory {\bf 130} (2010), no. 11, 2528--2541.

\bibitem{buc2} A. Bucur, C. David, B. Feigon, M. Lal\'in, \emph{Biased statistics for traces of cyclic $p$-fold covers over finite fields}, WIN—women in numbers, 121--143, Fields Inst. Commun., {\bf 60}, Amer. Math. Soc., Providence, RI, 2011.

\bibitem{buc3} A. Bucur, C. David, B. Feigon, M. Lal\'in, K. Sinha, \emph{Distribution of zeta zeroes of Artin-Schreier curves}, to appear in Math. Res. Lett.

\bibitem{bk} A. Bucur, K. Kedlaya, \emph{The probability that a complete intersection is smooth}, J. Theor. Nombres Bordeaux {\bf 24} (2012), no. 3, 541--556.


\bibitem{Del74} P. Deligne, \emph{La conjecture de Weil. I}, Inst. Hautes \'{E}tudes Sci. Publ. Math. {\bf 43} (1974), 273--307.

\bibitem{ent} A. Entin, \emph{On the distribution of zeroes of Artin-Schreier L-functions}, Geom. Funct. Anal. {\bf 22} (2012), 1322--1360.



\bibitem{fai} D. Faifman, Z. Rudnick, \emph{Statistics of the zeros of zeta functions in families of hyperelliptic curves over a finite field}, Compos. Math. {\bf 146} (2010), no. 1, 81--101.



\bibitem{kat} N. M. Katz, P. Sarnak, ``Random Matrices, Frobenius Eigenvalues, and Monodromy'', Amer. Math. Soc. Colloq. Publ., vol. 45, American Mathematical Socitey, Providence, RI, 1999.


\bibitem{kur} P. Kurlberg, Z. Rudnick, \emph{the fluctuations in the number of points on a hyperelliptic curve over a finite field}, J. Number Theory {\bf 129} (2009), no. 3, 580--587.

\bibitem{kur2} P. Kurlberg, I. Wigman, \emph{Gaussian point count statistics for families of curves over a fixed finite field}, Int. Math. Res. Not. IMRN 2011, no. 10, 2217--2229.

\bibitem{mon} H. L. Montgomery, ``Ten lectures on the interface between analytic number theory and harmonic analysis''. CBMS Regional Conference Series in Mathematics, {\bf 84}. American Mathematical Society, Providence, RI, 1994.

\bibitem{mor} C. Moreno, ``Algebraic curves over finite fields'', Cambridge Tracts in Mathematics {\bf 97}, Cambridge University Press, 1991.

\bibitem{neu} J. Neukirch, ``Algebraic Number Theory''. Translated from the 1992 German original and with a note by Norbert Schappacher. With a foreword by G. Harder. Grundlehren der Mathematischen Wissenschaften (Fundamental Principles of Mathematical Sciences), {\bf 322}. Springer-Verlag, Berlin, 1999.

\bibitem{ros} M. Rosen, ``Number theory in function fields''. Graduate Texts in Mathematics, {\bf 210}. Springer-Verlag, New York, 2002.









\bibitem{tay} M.J. Taylor, \emph{On the equidistribution of Frobenius in cyclic extensions of a number field}, J. London Math. Soc. {\bf 29} (1984), no. 2, 211--223.

\bibitem{wei} A. Weil, \emph{Sur les Courbes Alg\'ebriques et les Vari\'et\'es qui s'en D\'eduisent}, Publ. Inst. Math. Univ. Strasbourg {\bf 7} (1945),  Hermann et Cie., Paris, 1948.


\bibitem{woo1} M.M. Wood, \emph{The distribution of the number of points on trigonal curves over $\F$}, Int. Math. Res. Not. IMRN 2012, doi: 10.1093/imrn/rnr256.

\bibitem{woo} M.M. Wood, \emph{On the probabilities of local behaviors in abelian field extensions}, Compos. Math. {\bf 146} (2010), no. 1, 102--128.

\bibitem{wri} D. Wright, \emph{Distribution of discriminants of abelian extensions}, Proc. London Math. Soc. (3) {\bf 58} (1989), no. 1, 17--50.

\bibitem{xio} M. Xiong, \emph{Statistics of the zeros of zeta functions in a family of curves over a finite field}, Int. Math. Res. Not. IMRN 2010, no. 18, 3489--3518.
\end{thebibliography}
\end{document}